\documentclass[a4paper,12pt,twoside]{article}
\usepackage[latin1]{inputenc}
\usepackage{amsfonts,amssymb,amsmath,exscale}
\usepackage[dvips]{graphicx,color,psfrag}
\usepackage{fancyhdr}
\usepackage{rotating}
\usepackage{pifont}
\usepackage[hyperref]{xcolor}
\definecolor{myblue}{RGB}{0,0,139}
\usepackage[colorlinks=true,citecolor=blue,urlcolor=blue]{hyperref}
\hypersetup{bookmarksnumbered}
\usepackage[square]{natbib}
\usepackage{tensor}
\usepackage{setspace}
\usepackage{bbm}
\usepackage[font=small,labelfont=bf]{caption}
\fancypagestyle{plain}{%
  \fancyhead{}%
  \fancyfoot[c]{\sffamily\thepage}%
}
\makeatletter
\def\cleardoublepage{\clearpage\if@twoside \ifodd\c@page\else
  \hbox{}
  \vspace*{\fill}
  \thispagestyle{empty}
  \newpage
  \if@twocolumn\hbox{}\newpage\fi\fi\fi}
\usepackage{algorithm2e}

\usepackage{amsmath,amssymb}
\renewcommand{\citep}[1]{(\cite{#1}{})}

\DeclareRobustCommand{\BA}{{\boldsymbol{\mathnormal A}}}
\DeclareRobustCommand{\BB}{{\boldsymbol{\mathnormal B}}}
\DeclareRobustCommand{\BC}{{\boldsymbol{\mathnormal C}}}

\DeclareRobustCommand{\BF}{{\boldsymbol{\mathnormal F}}}

\DeclareRobustCommand{\BJ}{{\boldsymbol{\mathnormal J}}}
\DeclareRobustCommand{\BK}{{\boldsymbol{\mathnormal K}}}

\DeclareRobustCommand{\BP}{{\boldsymbol{\mathnormal P}}}

\DeclareRobustCommand{\BR}{{\boldsymbol{\mathnormal R}}}

\DeclareRobustCommand{\BX}{{\boldsymbol{\mathnormal X}}}

\DeclareRobustCommand{\Ba}{{\boldsymbol{\mathnormal a}}}

\DeclareRobustCommand{\Bd}{{\boldsymbol{\mathnormal d}}}

\DeclareRobustCommand{\Bf}{{\boldsymbol{\mathnormal f}}}

\DeclareRobustCommand{\Bn}{{\boldsymbol{\mathnormal n}}}

\DeclareRobustCommand{\Bt}{{\boldsymbol{\mathnormal t}}}

\DeclareRobustCommand{\Bx}{{\boldsymbol{\mathnormal x}}}

\DeclareRobustCommand{\calB}{{\mathcal B}}
\DeclareRobustCommand{\calD}{{\mathcal D}}
\DeclareRobustCommand{\calH}{{\mathcal H}}
\DeclareRobustCommand{\calL}{{\mathcal L}}
\DeclareRobustCommand{\calP}{{\mathcal P}}
\DeclareRobustCommand{\calR}{{\mathcal R}}
\DeclareRobustCommand{\calT}{{\mathcal T}}
\DeclareRobustCommand{\calW}{{\mathcal W}}

\DeclareMathAlphabet{\Inbb}{U}{bbmss}{m}{n}

\DeclareMathAlphabet{\gothic}{U}{euf}{m}{n}

\DeclareMathAlphabet{\Bgothic}{U}{euf}{b}{n}

\newcommand{\tr}{\mathop{\operator@font tr}}
\newcommand{\dev}{\mathop{\operator@font dev}}

\DeclareMathAlphabet{\Ibb}{U}{msb}{m}{n}

\newcommand{\een}{\]\@ignoretrue}

\DeclareRobustCommand{\Bvarphi}{{\boldsymbol{\varphi}}}
\DeclareRobustCommand{\Bomega}{{\boldsymbol{\omega}}}
\DeclareRobustCommand{\Bxi}{{\boldsymbol{\xi}}}
\renewcommand{\div}{\mathop{\operator@font div}}

\newcounter{SetOfEq}

\newcommand{\assemble}{\mathop{\textbf{\Large{\sf A}}}}

\begin{document}

\thispagestyle{empty}

\begin{center}
\textbf{\large Transient stability analysis of composite hydrogel structures \\[0.2cm]
based on a minimization-type variational formulation}

\medskip

S.~Sriram, E.~Polukhov \& M.-A.~Keip%
\footnote{Corresponding author: marc-andre.keip@mechbau.uni-stuttgart.de, Phone: +49 711 685 66233, Fax: +49 711 685 66347}

\medskip

Institute of Applied Mechanics \\
Department of Civil and Environmental Engineering \\
University of Stuttgart, Stuttgart, Germany
\end{center}
\medskip

\noindent \textbf{Abstract.} We employ a canonical variational framework for the predictive characterization of structural instabilities that develop during the diffusion-driven transient swelling of hydrogels under geometrical constraints. The variational formulation of finite elasticity coupled with Fickian diffusion has a \textit{two-field minimization} structure, wherein the deformation map and the fluid-volume flux are obtained as minimizers of a time-discrete potential involving internal and external energetic contributions. Following spatial discretization, the minimization principle is implemented using a conforming Q\textsubscript{1}RT\textsubscript{0} finite-element design, making use of the lowest-order Raviart--Thomas-type interpolations for the fluid-volume flux. To analyze the structural stability of a certain equilibrium state of the gel satisfying the minimization principle, we apply the \textit{local stability criterion} on the incremental  potential, which is based on the idea that a \textit{stable} equilibrium state  has the lowest potential energy among all possible states within an infinitesimal neighborhood. Using this criterion, it is understood that bifurcation-type structural instabilities are activated when the coupled global finite-element stiffness matrix \textit{loses its positive definiteness}. This concept is then applied to determine the onset and nature of wrinkling instabilities occurring in a pair of representative film-substrate hydrogel systems. In particular, we analyze the dependencies of the critical buckling load and mode shape on the system geometry and material parameters.


\pagestyle{fancy}
\pagenumbering{arabic}

\section{Introduction}

In view of the recent surge in demand for soft multi-functional materials, hydrogels have gained a lot of importance owing to their increasing scientific and industrial applications. Essentially, they are soft hydrophilic elastomers that exhibit a great degree of swelling upon absorbing a diffusing fluid. In the swollen state, the presence of cross-links between the elastomer chains provides them with favourable mechanical properties. This together with their biocompatibility make hydrogels highly sought after for several biomedical applications such as drug delivery systems and contact lenses, refer \citet{calo2015}.       

When the free swelling of hydrogels is arrested by suitable mechanical constraints, a variety of structural instability patterns are observed. Experimental studies by \citet{tanaka1987}, \citet{trujillo2008} and \citet{guvendiren2010} have shown that creasing patterns begin to develop on polymer gel films attached to rigid surfaces when a certain critical linear expansion ratio is reached upon swelling. For a detailed review of more experimental studies on pattern formations in gels, we refer to \citet{dervaux2011}. More recently, \citet{liaw2019} have elucidated the various mechanisms through which surface instability patterns are activated in confined hydrogels. Advancements in the field of controlled pattern generation in hydrogel systems that respond to a wide range of external stimuli such as changes in solvent concentration, temperature, electric field etc.\ have also been reviewed by the authors. Through predictive modeling of such mechanical instabilities, it is possible to design hydrogel systems over a wide range of material parameters and geometries that are tuned to buckle and produce selected surface morphologies, which can then be exploited for specific applications. In this regard, we refer to the works of \citet{yin2009}, \citet{chen2010}, \citet{yang2010} and \citet{huang2013}.     

Over the past decade, there has been a considerable amount of research focused on developing analytical/semi-analytical methods for the stability problem of soft elastomers subjected to constrained growth. In the works of \citet{cao2011,cao2012} and \citet{jin2015}, an equilibrium bifurcation analysis has been adopted to study the onset of sinusoidal wrinkles on the surfaces of neo-Hookean film-substrate bilayers under plane-strain compression over a broad range of modulus and thickness ratios. In the limit when both layers have the same modulus, wrinkles are found to be highly unstable and undergo a dynamic transition to creases. A similar analysis on elastic cylindrical bilayers has been conducted by \citet{li2011} and \citet{xie2014}. Here, a multiplicative split of the deformation gradient into separate elastic deformation and growth tensors has been considered, following which a linear perturbation analysis of the incremental equilibrium equation is carried out. The latter have proven that internal pressure and surface tension have a stabilizing effect on the cylindrical bilayer. For thin films growing on the outer surfaces of soft cylindrical substrates, \citet{jin2018} and \citet{jia2018} have presented asymptotic solutions predicting the onset and nature of wrinkling instabilities for a broad range of film-substrate modulus ratios. Such asymptotic solutions, which are derived using the concept of order analysis, have been shown to predict the buckling characteristics with good accuracy, particularly when the film and substrate layers have comparable shear moduli. Comparing the results for the cylindrical geometry to that of flat bilayers with similar thicknesses, material properties and boundary conditions, the latter have shown that curvature tends to have a stabilizing effect on cylindrical bilayers. \citet{jin2019} have extended the concept of asymptotic solutions to study the dependence of wrinkle amplitude on the geometrical and material parameters of elastic cylindrical bilayers. For an analytical investigation concerning the bifurcation characteristics of growing spherical shells under external pressure, we refer to \citet{benamar2005}. 

Common to the above works is the consideration of a purely elastic system with no chemical coupling due to species diffusion. Semi-analytical investigations on the diffusion-driven swelling-induced instabilities occurring in constrained hydrogels have been carried out by \citet{kang2010}, \citet{xiao2012} and \citet{wu2013,wu2017onset} for different shapes and material profiles. The mechanism that drives the formation of wrinkles and creases in such constrained hydrogels is comparable to that which activates instabilities in equivalent hyperelastic systems having similar boundary conditions. As a result, one could expect similar trends for the critical buckling characteristics in both cases with respect to specimen geometries and material parameters.

One of the major disadvantages associated with analytical modeling of instability phenomena is that its applicability is restricted to simple load cases and system geometries. In contrast, computational methods provide the possibility to explore the bifurcation characteristics of hydrogel systems with more complex shapes, material inhomogeneities and loading conditions. An important prerequisite for this is that the proposed computational model should accurately and efficiently capture the finite-strain swelling response of the hydrogel under fluid diffusion. Using a combination of neo-Hookean and Flory--Rehner-type energy functions for modeling the coupled response, \citet{hong2009,hong2008}, \citet{bouklas2015} and \citet{chester2015} have proposed direct methods based on the principle of virtual work to model the diffusion-driven swelling of hydrogel structures. \citet{ilseng2019} have extended the above modeling framework by carrying out a stability analysis to study the formation of wrinkles on the surfaces of laterally-constrained layered hydrogel plates. In addition to investigating the critical swelling ratios and critical wrinkle wavelengths for various initial geometries of the composite plate, the influence of diffusion coefficient has also been studied. 

In the context of multi-field problems, such direct methods often result in unsymmetric system matrices (as seen in \citet{bouklas2015}) upon finite-element discretization, which hamper the computational efficiency of the model. An alternate strategy for modeling multi-physics problems is to adopt a \textit{variational-based} approach, wherein the solution of the coupled boundary value problem is obtained by optimizing a certain objective functional. In a finite-element context, the global stiffness matrix is obtained as the second derivative of the discretized objective functional and is therefore \textit{inherently symmetric}. We refer to the works of \citet{miehe2014}, \citet{boeger2017} and \citet{teichtmeister2019}, where a variety of variational minimization and saddle-point formulations for standard Fickian-type and gradient-extended Cahn--Hilliard-type species diffusion in elastic solids have been developed and compared. Such formulations have been shown to efficiently capture phenomena such as diffusion-induced finite deformation of hydrogel rods and pressure-induced diffusion in indentation-type problems. Recently, \citet{zheng2020} have modeled the transient coupled response of hydrogels by adopting a mixed isogemetric analysis approach. Problems such as the free swelling of cylindrical hydrogel structures immersed in a solvent and the transient swelling of hydrogel blocks subjected to local solvent injection have been successfully studied by discretizing the geometries with mixed isogeometric elements. Such elements make use of Non-Uniform Rational B-Spline (NURBS) basis functions to circumvent volumetric locking effects caused by the incompressibility of polymer chains and solvent molecules.  

The procedure to incorporate a structural stability analysis into the framework of a variational formulation has been outlined by \citet{geymonat1993}, \citet{miehe2002,mieheval2014}, \citet{schroder2017} and \citet{dortdivanlioglu2018} for both single- and multi-field problems. In the work of the latter, diffusion-driven structural instabilities occurring in constrained hydrogels have been modeled using a saddle-point formulation consisting of the deformation and chemical-potential fields as independent variables. The saddle-point nature of the problem leads to an \textit{indefinite} stiffness matrix, which makes the accompanying stability analysis relatively less straightforward. In this work, we adopt a minimization-based formulation, having the deformation and fluid-volume flux as primary fields, since it provides a convenient platform to carry out a structural stability analysis, owing to the inherent \textit{positive-definiteness} of the global stiffness matrix for a stable chemo-mechanical equilibrium state. We refer to \citet{polukhov2020} and the references cited therein for an overview of the advantages associated with a minimization formulation over a saddle-point approach. For finite elasticity coupled with Fickian-type diffusion, the minimization-based approach demands a more complex and non-trivial spatial discretization scheme for the volume-flux field due to $\calH(\text{Div})$-conformity requirement. This will be discussed in greater detail in Appendix~\ref{app:A}.

The present work is structured as follows. In Section 2, we begin by presenting the strong form of the coupled initial boundary value problem, i.e.\ the governing kinematic, equilibrium and constitutive equations together with the boundary conditions. A rate-type potential functional, formulated in terms of the deformation rate and volume-flux fields, is postulated. Following an implicit Euler time integration, the resulting incremental potential is spatially discretized using conforming Q\textsubscript{1}RT\textsubscript{0} finite elements, which use the lowest-order vectorial Raviart--Thomas shape functions for interpolating the nodal fluxes (cf.~Appendix A). A local stability criterion is then introduced to judge the structural stability of a certain chemo-mechanical equilibrium state based on an eigenvalue analysis of the coupled global stiffness matrix. Following a brief description of the adopted constitutive model for hydrogels, we implement the discrete variational framework to investigate the critical buckling characteristics of a pair of representative composite hydrogel structures in Section 3. In particular, variations in the critical conditions and corresponding buckling patterns with respect to geometrical and material parameters are studied. Finally, a brief summary of the work is presented in Section 4.      

\section{Variational stability analysis in the framework of a finite-strain chemo-mechanical minimization formulation}

We begin by presenting the continuum-mechanical closure problem for finite-strain chemo-mechanics. Subsequently, a rate-type variational formulation of the coupled problem is discussed in a continuous setting, following which discretizations in time and space are carried out. The following presentation is based on the works of \citet{miehe2014} and \citet{boeger2017,boger2017minimization}.

\subsection{Primary fields and equilibrium equations for the coupled problem}

In the present context, diffusion-driven finite elasticity is modeled as a two-field problem consisting of the \textit{deformation map} $\Bvarphi: \calB_0 \times \calT \to \calB_t \subset \calR^3$, mapping a material point $\BX \in \calB_0$ onto its counterpart $\Bx \in \calB_t$ in the deformed configuration, and \textit{fluid-volume flux} $\mathbb{H}: \calB_0 \times \calT \to \calR^3$ as global primary fields. The deformation gradient 
\begin{equation}
	\label{eq:2.1}
	\BF:=\nabla_\BX \Bvarphi (\BX,t)\quad \text{with}\quad J:=\text{det}[\BF]>0	
	\end{equation}
represents the fundamental measure of local mechanical stretch in a geometrically nonlinear setting. The fluid-volume flux $\mathbb{H}(\BX,t)$ represents the volume of fluid diffusing across a unit area element $\mathrm{d}\BA$ of the undeformed continuum per unit time. It is a pull-back of the Eulerian quantity $\mathbbm{h}(\Bx,t)$ representing the fluid flux across an area element $\mathrm{d}\Ba$ of the deformed configuration, i.e.
\begin{equation}
	\label{eq:2.2}
	\mathbb{H}\cdot\mathrm{d}\BA=\mathbbm{h}\cdot\mathrm{d}\Ba \,\,\implies\,\,\mathbb{H}:=J\BF^{-1}\mathbbm{h}\,.
\end{equation}
Additionally, we have the fluid-concentration field $s:\calB_0 \times \calT \to \calR^+$ representing locally the volume of fluid in an infinitesimal Lagrangian volume-element relative to its unswollen volume. It is related to the primary field $\mathbb{H}$ by the local balance of fluid volume
\begin{equation}
	\label{eq:2.3}
	\dot{s}=-\text{Div}[\mathbb{H}]\,,
\end{equation} 
where $\dot{(\cdot)}$ represents the material time derivative of the quantity $(\cdot)$. 

In the absence of mechanical body forces and inertial effects, we have the balance equations for linear and angular momenta given by
\begin{equation}
	\label{eq:2.4}
	\text{Div}[\BP]=\textbf{0} \quad\text{and}\quad \text{skew}[\BP\BF^T]=\textbf{0}\,,
\end{equation}
respectively. In the above equations, $\BP$, the first Piola-Kirchhoff stress tensor, is the energetically conjugate variable to the deformation gradient $\BF$. The set of equations (\ref{eq:2.3}) and (\ref{eq:2.4}) represents the equilibrium equations governing the coupled diffusion-deformation process. For closure of the chemo-mechanical problem, we need constitutive laws describing the coupled material response, which are to be formulated in such a way that the principle of thermodynamic irreversibility is satisfied \textit{a priori}. Referring to \citet{boeger2017} and \citet{Miehe2015}, the dissipation postulate is split into separate restrictions on the local and convective parts that are to be individually satisfied, i.e.
\begin{equation}
	\label{eq:2.5}
	\calD_{loc}:=\BP:\dot{\BF}+\mu\dot{s}-\dot{\psi}\geq 0\quad \text{and} \quad
	\calD_{dif}:=-\mathbb{H}\cdot\nabla\mu \geq 0\,.
\end{equation}   
\begin{figure}[t]
	\centering
	\psfrag{n0}[c][c]{$\Bn_0$}
	\psfrag{XinB}[c][c]{$\BX \in \calB_0$}
	\psfrag{phi}[l][l]{$\Bvarphi$}
	\psfrag{h}[l][l]{$\mathbb{H}$}
	\psfrag{tbar}[l][l]{$\BP\,\Bn_0=\bar{\Bt}_0$}
	\psfrag{hbar}[l][l]{$\mathbb{H}\cdot\Bn_0=\bar{h}_0$}
	\psfrag{phibar}[c][c]{$\Bvarphi=\bar{\Bvarphi}$}
	\psfrag{mubar}[c][c]{$\mu=\bar{\mu}$}
	\psfrag{dbphi}[c][c]{$\partial\calB_0^{\,\Bvarphi}$}
	\psfrag{dbmu}[c][c]{$\partial\calB_0^{\,\mu}$}
	\psfrag{dbh}[c][c]{$\partial\calB_0^{\,h}$}
	\psfrag{dbt}[c][c]{$\partial\calB_0^{\,\Bt}$}
	\includegraphics[width=0.75\textwidth]{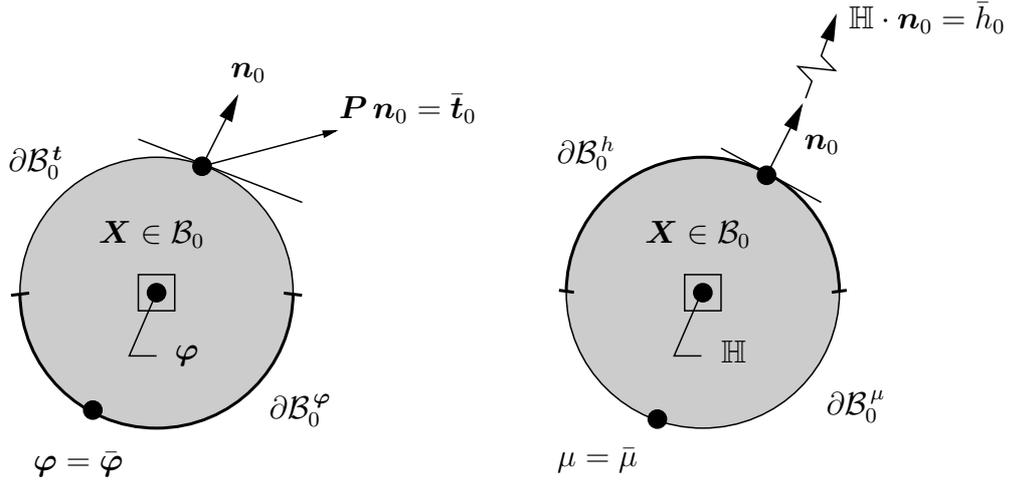}
	\captionsetup{width=0.9\textwidth}
	\caption{\footnotesize \textit{Primary fields of the coupled problem and the boundary conditions.} Within the minimization-based variational framework, finite elasticity coupled with diffusion is modeled using two primary fields, namely the deformation $\Bvarphi$ and the fluid-volume flux $\mathbb{H}$. The total boundary $\partial\calB_0$ is split into $\partial\calB_0^{\,\Bvarphi}$ and $\partial\calB_0^{\,\Bt}$ for the mechanical fields and $\partial\calB_0^{\,h}$ and $\partial\calB_0^{\,\mu}$ for the chemical fields such that $\partial\calB_0^{\,\Bvarphi}\cap\partial\calB_0^{\,\Bt}=\partial\calB_0^{\,h}\cap\partial\calB_0^{\,\mu}=\varnothing$. The vector $\Bn_0(\BX)$ represents the unit outward normal to the boundary $\partial\calB_0$ of the continuum in the reference configuration. Conditions for the various fields on their respective boundaries are shown.}
	\label{fig:fields}
\end{figure}
Here, $\mu(\BX,t)$ represents the chemical potential and is the energetically conjugate variable to the fluid concentration $s$. Following a local theory of grade one, the former restriction leads us to the constitutive equations
\begin{equation}
	\label{eq:2.6}
	\BP=\partial_{\BF}\hat{\psi}(\BF,s) \quad \text{and} \quad 
	\mu = \partial_{s}\hat{\psi}(\BF,s)
\end{equation}
for an \textit{objective} free-energy function $\hat{\psi}(\BF,s)$. The latter restriction is fulfilled by introducing a dissipation potential $\hat{\phi}(\mathbb{H};\BF,s)$ at a known chemo-mechanical state $\left\lbrace  \BF,s\right\rbrace$ such that we have the constitutive equation
\begin{equation}
	\label{eq:2.7}
	\nabla\mu = -\partial_{\,\mathbb{H}}\hat{\phi}(\mathbb{H};\BF,s)\,.
\end{equation}
Plugging (\ref{eq:2.7}) into (\ref{eq:2.5})\textsubscript{2}, we observe that the inequality is automatically satisfied if the dissipation potential $\hat{\phi}$ is formulated as a \textit{convex} homogeneous function in $\mathbb{H}$ such that $\partial_{\,\mathbb{H}}\hat{\phi}(\mathbb{H}=\textbf{0})=\textbf{0}$. Convexity of $\hat{\phi}$ is a sufficient condition for the minimization-based structure of the two-field variational formulation discussed subsequently. Finally, to complete the definition of the coupled initial boundary value problem, we prescribe conditions for the primary fields and their dual variables on the mechanical and chemical parts of the boundary. For a well-posed problem, the split in boundary $\partial\calB_0$ of the reference continuum should satisfy the conditions
$\partial\calB_0^{\,\Bvarphi}\cup\partial\calB_0^{\,\Bt}=\partial\calB_0=\partial\calB_0^{\,h}\cup\partial\calB_0^{\,\mu}$ and
$\partial\calB_0^{\,\Bvarphi}\cap\partial\calB_0^{\,\Bt}=\varnothing=\partial\calB_0^{\,h}\cap\partial\calB_0^{\,\mu}$ on the mechanical and chemical parts, as shown in Fig.~\ref{fig:fields}. We then have
\begin{equation}
	\label{eq:2.9}
	\Bvarphi = \bar{\Bvarphi}\,\,\,\text{on}\,\,\,\partial\calB_0^{\,\Bvarphi}\quad\text{and}\quad
	\mathbb{H}\cdot\Bn_0= \bar{h}_0\,\,\,\text{on}\,\,\,\partial\calB_0^{\,h} 
\end{equation} 
as the Dirichlet boundary conditions and 
\begin{equation}
	\label{eq:2.10}
	\BP\,\Bn_0=\bar{\Bt}_0\,\,\,\text{on}\,\,\,\partial\calB_0^{\,\Bt}\quad\text{and}\quad  \mu=\bar{\mu}\,\,\,\text{on}\,\,\,\partial\calB_0^{\,\mu}
\end{equation}
as the Neumann boundary conditions within the minimization framework. In addition to (\ref{eq:2.9}) and (\ref{eq:2.10}), we also have the initial condition 
$s(\BX,0)=s_0 \,\,\, \text{in}\,\,\,\calB_0$
accompanying the temporal differential equation (\ref{eq:2.3}) for the fluid concentration.

\subsection{A continuous rate-type minimization principle for the evolution problem}

Diffusion being a transient phenomenon, a rate-type potential functional governing the coupled problem is constructed as the difference between the rates of internal and external potentials, i.e.
\begin{equation}
	\label{eq:2.12}
	\dot{\Pi}(\dot{\Bvarphi},\mathbb{H})=\dot{\Pi}_{int}(\dot{\BF},\text{Div}[\mathbb{H}],\mathbb{H}) - \dot{\Pi}_{ext}(\dot{\Bvarphi},\mathbb{H})\,.
\end{equation} 
The rate of internal potential $\dot{\Pi}_{int}$ constitutes two functional contributions, namely the rate of energy storage and the rate of dissipation functionals, and is expressed as
\begin{equation}
	\label{eq:2.13}
	\dot{\Pi}_{int}(\dot{\BF},\text{Div}[\mathbb{H}],\mathbb{H}) =
	\frac{\mathrm{d}}{\mathrm{d}t}\int_{\calB_0}\hat{\psi}(\BF,s)\,\mathrm{d}V + \int_{\calB_0}\hat{\phi}(\mathbb{H};\BF,s)\,\mathrm{d}V
\end{equation}
The balance of fluid volume (\ref{eq:2.3}) is incorporated as a constraint in (\ref{eq:2.13}). In the absence of mechanical body forces, the rate of external potential $\dot{\Pi}_{ext}$ consists of a mechanical power due to prescribed tractions $\bar{\Bt}_0$ on $\partial\calB_0^{\,\Bt}$ and a power due to fluid diffusion driven by a prescribed chemical potential $\bar{\mu}$ on $\partial\calB_0^{\,\mu}$
\begin{equation}
	\label{eq:2.14}
	\dot{\Pi}_{ext}(\dot{\Bvarphi},\mathbb{H})=\dot{\Pi}_{ext}(\dot{\Bvarphi})+\dot{\Pi}_{ext}(\mathbb{H})
        =\int_{\partial\calB_0^{\,\Bt}}\bar{\Bt}_0\cdot\dot{\Bvarphi}\,\,\mathrm{d}A-\int_{\partial\calB_0^{\,\mu}}\bar{\mu}\,\mathbb{H}\cdot\Bn_0\,\mathrm{d}A\,. 
\end{equation}
The negative sign in the term $\dot{\Pi}_{ext}(\mathbb{H})$ is due to the fact that $\mathbb{H}\cdot\Bn_0$ by itself represents an outflux of fluid molecules. Hence, for the case where fluid diffuses \textit{into} the solid continuum, the integrand turns out to be negative, resulting in a positive fluid power term. In the representative numerical examples that will be discussed in subsequent sections, the deforming surfaces of the specimen under study are considered to be \textit{traction-free} at every instant of time. As a result, the overall rate of potential (\ref{eq:2.12}) reduces to the form
\begin{equation}
	\label{eq:2.15}
	\dot{\Pi}(\dot{\Bvarphi},\mathbb{H})= \int_{\calB_0}\left[ \partial_{\BF}\hat{\psi}:\dot{\BF}-\partial_{s}\hat{\psi}\,\text{Div}[\mathbb{H}]+\hat{\phi}(\mathbb{H};\BF,s)\right] \mathrm{d}V + \int_{\partial\calB_0^{\,\mu}}\bar{\mu}\,\mathbb{H}\cdot\Bn_0\,\mathrm{d}A\,.
\end{equation} 
In the works of \citet{boeger2017} and \citet{teichtmeister2019}, it has been shown that the governing equilibrium equations (\ref{eq:2.3}) and (\ref{eq:2.4}) for finite elasticity coupled with Fickian diffusion and the prescribed Neumann boundary conditions (\ref{eq:2.10}) are recovered as Euler-Lagrange equations of the rate-type variational formulation (\ref{eq:2.15}), proving that the proposed variational framework is consistent with the underlying physics of the coupled diffusion-deformation process. For a convex dissipation potential, the primary fields of the continuous formulation are obtained by the minimization principle (\cite{miehe2014}, \cite{boeger2017}, \cite{teichtmeister2019}) 
\begin{equation}
	\label{eq:2.16}
\boxed{\left\lbrace\dot{\Bvarphi}^*,\mathbb{H}^* \right\rbrace = \text{Arg}\left\lbrace \inf_{\dot{\Bvarphi}}\inf_{\mathbb{H}}\dot{\Pi}(\dot{\Bvarphi},\mathbb{H})\right\rbrace}
\end{equation}  
subject to the conditions $\dot{\Bvarphi}=\bar{\dot{\Bvarphi}}$ on $\partial\calB_0^{\,\Bvarphi}$ and $\mathbb{H}\cdot\Bn_0=\bar{h}_0$ on $\partial\calB_0^{\,h}$.

\subsection{Space-time-discrete form of the minimization principle}

\textit{Temporal discretization of the rate-type potential.}
In order to recover the primary fields $\Bvarphi$ and $\mathbb{H}$, we apply a fully-implicit backward Euler integration scheme to the rate-type potential described in (\ref{eq:2.15}) and to the evolution equation (\ref{eq:2.3}) for the fluid concentration over a time step $\tau:=t_{n+1}-t_n$. It is assumed that the minimization principle holds in this chosen discrete time step. This results in an incremental form of the coupled potential given by
\begin{equation}
	\label{eq:2.17}
	\begin{aligned}
	\Pi^\tau(\Bvarphi,\mathbb{H})&=\int_{\calB_0}\pi^\tau(\BF,\text{Div}[\mathbb{H}],\mathbb{H})\,\mathrm{d}V +
	\int_{\partial\calB_0^{\,\mu}}\tau\bar{\mu}\,\mathbb{H}\cdot\Bn_0\,\mathrm{d}A\,,\\[1ex]
	\hspace{-21mm}\text{with}\quad
	\pi^\tau(\BF,\text{Div}[\mathbb{H}],\mathbb{H})&=\hat{\psi}(\BF,s_n-\tau\text{Div}[\mathbb{H}])-\hat{\psi}(\BF_n,s_n)+\tau\hat{\phi}(\mathbb{H};\BF_n,s_n)\,.
	\end{aligned}
\end{equation}
Quantities without a subscript represent values at the current time $t_{n+1}$. The incremental dissipation potential $\tau\hat{\phi}(\mathbb{H};\BF_n,s_n)$ in (\ref{eq:2.17})\textsubscript{2} is a function of the current volume flux $\mathbb{H}$ and is evaluated at a known state of deformation and swelling $\left\lbrace \BF_n,s_n\right\rbrace $ in order to ensure variational consistency. The primary fields $\Bvarphi$ and $\mathbb{H}$ at the current time $t_{n+1}$ are then the minimizers of the incremental potential $\Pi^\tau$, i.e.
\begin{equation}
	\label{eq:2.18}
\boxed{\left\lbrace \Bvarphi^*,\mathbb{H}^*\right\rbrace= \text{Arg}
	\left\lbrace\inf_{\Bvarphi\in\calW_\Bvarphi} \inf_{\mathbb{H}\in\calW_\mathbb{H}}\Pi^\tau(\Bvarphi,\mathbb{H})\right\rbrace} \,, 
\end{equation}
where the set of admissible deformations $\calW_\Bvarphi$ and fluid-volume fluxes $\calW_\mathbb{H}$ read
$\calW_\Bvarphi:=\left\lbrace\Bvarphi \in \calH^1(\calB_0) \,|\, \Bvarphi=\bar{\Bvarphi}\,\,\text{on}\,\,\partial\calB_0^{\,\Bvarphi} \right\rbrace$
and
$\calW_\mathbb{H}:=\left\lbrace\mathbb{H} \in \calH(\text{Div},\calB_0)\, |\, \mathbb{H}\cdot\Bn_0=\bar{h}_0\,\,\text{on}\,\,\partial\calB_0^{\,h} \right\rbrace$.

The incremental internal potential density $\pi^\tau$ introduced in (\ref{eq:2.17}) is a function of the current state of the continuum represented by the constitutive-state array $\mathfrak{C}:=[\BF,\text{Div}[\mathbb{H}],\mathbb{H}]^T$. The first and second derivatives of $\pi^\tau$ with respect to $\mathfrak{C}$ yield the array of driving forces $\mathfrak{D}$ and tangent moduli $\mathbb{C}$ given by
\begin{equation}
	\label{eq:2.20}
	\mathfrak{D}:=\partial_\mathfrak{C}\pi^\tau=
	\begin{bmatrix}
	\,\partial_{\BF}\hat{\psi} \\ -\tau\partial_{s}\hat{\psi} \\ \tau\partial_{\,\mathbb{H}}\hat{\phi}\,
	\end{bmatrix}
	\quad \text{and} \quad
	\mathbb{C}:=\partial^2_{\mathfrak{C}\mathfrak{C}}\pi^\tau=
	\begin{bmatrix}
	\partial_{\BF \BF}^2 \hat{\psi}\quad& -\tau\partial_{\BF s}^2 \hat{\psi}\quad& \cdot \\
	-\tau\partial_{s \BF}^2 \hat{\psi}\quad& \tau^2 \partial_{s s}^2 \hat{\psi}\quad& \cdot \\
	\cdot \quad& \cdot \quad& \tau \partial_{\mathbb{H}\mathbb{H}}^2 \hat{\phi}
	\end{bmatrix}\,,
\end{equation}
respectively. The arrays $\mathfrak{D}$ and $\mathbb{C}$ enter the finite-element implementation of the space-time-discrete minimization principle discussed below. 

\textit{Spatial discretization in a finite-element setting.}
In the following, we discuss the generalized procedure for the finite-element implementation of the minimization-based variational framework for two-dimensional problems. Consider a spatial discretization of the reference body $\calB_0$ into $n$ discrete two-dimensional finite elements resulting in a discretized configuration 
$\calB_0^h:=\bigcup_{e=1}^n\, \calB_0^e$.
Let $\Bd^{\,e}:=\left\lbrace\Bd^{\,\Bvarphi},\Bd^{\,h} \right\rbrace $ denote the vector containing the displacement and flux degrees of freedom for a finite element $e$. Then, the constitutive-state array $\mathfrak{C}$ for this element reads
$\mathfrak{C}=\BB^e(\BX)\,\Bd^{\,e}$,
where $\BB^e$ represents an element-level matrix mapping the nodal degrees of freedom to the constitutive state. It is a function of the set of material points $\BX$ of the body that are covered by the element~$e$.

The space-discrete form of the incremental potential $\Pi^\tau$ (\ref{eq:2.17})\textsubscript{1} is then given by
\begin{equation}
	\label{eq:2.23}
	\Pi^{\tau,h}(\Bd)=\assemble_{e=1}^n\left\lbrace \int_{\calB_0^e} \pi^{\tau,h}(\BB^e\Bd^{\,e})\,\mathrm{d}V^e - \Pi_{ext}^{\tau,h}(\Bd^{\,e})\right\rbrace\,,
\end{equation}
where $\Bd$ is the global vector containing all the nodal degrees of freedom and $\assemble$ represents the finite-element assembly operator. The discrete incremental external potential $\Pi_{ext}^{\tau,h}$ results from the prescribed chemical potential at the boundary $\partial\calB_0^{\,\mu,e}$ of the finite element $e$. The solution vector $\Bd^*$ containing the current values of nodal displacements and fluxes is obtained by the space-time-discrete minimization principle
\begin{equation}
	\label{eq:2.24}
	\boxed{\Bd^*=\text{Arg}\left\lbrace\inf_\Bd \Pi^{\tau,h}(\Bd)\right\rbrace}\,. 
\end{equation} 
The necessary condition for the minimization principle (\ref{eq:2.24}) reads
$\BR(\Bd^*):=\Pi^{\tau,h}_{,\Bd}\,\big|_{\Bd=\Bd^*} \stackrel{!}{=} \textbf{0}$,	  
which upon simplification yields the equilibrium condition
\begin{equation}
	\label{eq:2.26}
	\assemble_{e=1}^n\left\lbrace\int_{\calB_0^e} (\BB^e)^T\mathfrak{D}^h\,\mathrm{d}V^e - \Bf_{ext}^e \right\rbrace \stackrel{!}{=} \textbf{0} .
\end{equation} 
The nonlinear equation (\ref{eq:2.26}) is solved by means of a global Newton--Raphson iterative scheme. The update for the global solution vector is expressed as
\begin{equation}
	\label{eq:2.27}
	\Bd^{\,i+1}=\Bd^{\,i}-(\BK^{\,i})^{-1}\BR^{\,i} \quad \text{with} \quad
	\BK(\Bd):=\Pi^{\tau,h}_{,\Bd\Bd}=\assemble_{e=1}^n\int_{\calB_0^e}(\BB^e)^T\mathbb{C}^h\BB^e\,\mathrm{d}V^e
\end{equation}
and is performed until $||\BR||$ is less than a certain tolerance. The coupled global stiffness matrix $\BK(\Bd)$ is \textit{symmetric} due to the fact that it is the second derivative of the discrete incremental potential $\Pi^{\tau,h}$. Additionally, since the global solution vector is obtained as a minimizer of $\Pi^{\tau,h}$, the matrix $\BK$ is inherently \textit{positive definite} for a stable equilibrium state $\Bd^*$.

\subsection{Variational-based structural stability analysis}
Consider a chemo-mechanical state $\Bomega_1^*:=\left\lbrace \Bvarphi_1^*, \mathbb{H}_1^*\right\rbrace $ satisfying the incremental two-field minimization principle (\ref{eq:2.18}) for admissible $\Bvarphi$ and $\mathbb{H}$. The state $\Bomega_1^*$ is considered to be \textit{globally stable} if the inequality
\begin{equation}
\label{eq:2.33}
\Pi^{\tau}(\Bomega_2)-\Pi^{\tau}(\Bomega_1^*) > 0
\end{equation}
holds for any other admissible state $\Bomega_2:=\left\lbrace\Bvarphi_2, \mathbb{H}_2 \right\rbrace $ satisfying (\ref{eq:2.18}). 

From the works of \citet{hill1957}, \citet{ball1976}, \citet{geymonat1993} and \citet{miehe2002}, it is seen that the global stability criterion given by (\ref{eq:2.33}) leads to the condition demanding \textit{strict convexity} of the free-energy function $\hat{\psi}$. The former has shown that a strictly convex energy function implies \textit{uniqueness of solution} of the associated boundary value problem, which is physically unacceptable in a finite-strain setting as it dismisses \textit{a priori} the possibility of bifurcation-type structural instabilities. 

We therefore consider the state $\Bomega_2$ to lie within an infinitesimal neighbourhood of $\Bomega_1^*$, i.e. $\Bomega_2:=\Bomega_1^* + \epsilon\,\delta\Bomega,$ giving rise to the condition
\begin{equation}
\label{eq:2.34}
\Pi^{\tau}(\Bomega_1^* + \epsilon\,\delta\Bomega)-\Pi^{\tau}(\Bomega_1^*) > 0
\end{equation} 
in order for $\Bomega_1^*$ to be a \textit{locally stable} state, refer \citet{mieheval2014}. Performing a Taylor-series expansion of $\Pi^\tau(\Bomega_2)$ about the state $\Bomega_1^*$ upto the second order term, (\ref{eq:2.34}) reduces to the form
\begin{equation}
\label{eq:2.35}
\frac{\mathrm{d}}{\mathrm{d}\epsilon}\bigg|_{\epsilon=0}\Pi^{\tau}(\Bomega_1^*+\epsilon\,\delta\Bomega)+\frac{1}{2!}\frac{\mathrm{d}^2}{\mathrm{d}\epsilon^2}\bigg|_{\epsilon=0}\Pi^{\tau}(\Bomega_1^*+\epsilon\,\delta\Bomega)+\cdot\cdot\cdot > 0\,,
\end{equation}
where the first and second directional derivatives in (\ref{eq:2.35}) represent the first and second variations of $\Pi^\tau$ at the state $\Bomega_1^*$, respectively. As a necessary condition of the minimization principle (\ref{eq:2.18}), the former term vanishes at the equilibrium state $\Bomega_1^*$. As a result, we arrive at the condition that the second variation of $\Pi^\tau$ at $\Bomega_1^*$ has to be positive, i.e. 
\begin{equation}
\label{eq:2.36}
\int_{\calB_0} 
\begin{bmatrix}
\delta \BF \\ \delta\,\text{Div}\left[\mathbb{H} \right] \\ \delta\mathbb{H} 
\end{bmatrix}
\cdot \begin{bmatrix}
\partial_{\BF \BF}^2 \hat{\psi}\quad& -\tau\partial_{\BF s}^2 \hat{\psi}\quad& \cdot \\
-\tau\partial_{s \BF}^2 \hat{\psi}\quad& \tau^2 \partial_{s s}^2 \hat{\psi}\quad& \cdot \\
\cdot \quad& \cdot \quad& \tau \partial_{\mathbb{H} \mathbb{H}}^2 \hat{\phi} 
\end{bmatrix}
\cdot \begin{bmatrix}
\delta \BF \\ \delta\,\text{Div}\left[\mathbb{H} \right] \\ \delta\mathbb{H} 
\end{bmatrix}
\, \mathrm{d}V\,\,\mathlarger{\mathlarger{\mathlarger{\bigg|}}}_{\Bomega_1^*} > 0 \,\,\text{.}
\end{equation} 

In a space-discrete setting, the above inequality can be rewritten using the finite-element arrays introduced previously in the form
\begin{equation}
\label{eq:2.37}
\assemble_{e=1}^n \int_{\calB_0^e} \delta \mathfrak{C}^T\,\mathbb{C}^h\,\delta \mathfrak{C} \,\,\mathrm{d}V^e\,\bigg|_{\Bd_1^*} > 0\,.
\end{equation}
Using the definition for the constitutive-state array, (\ref{eq:2.37}) simplifies to the statement
\begin{equation}
\label{eq:2.38}
\delta\Bd^T \assemble_{e=1}^n \int_{\calB_0^e} (\BB^e)^T\mathbb{C}^h \BB^e \,\mathrm{d}V^e \,\bigg|_{\Bd_1^*} \delta\Bd > 0 \quad \implies \quad \delta\Bd^T\, \BK(\Bd_1^*)\,\delta\Bd > 0 \,.
\end{equation}
In other words, the coupled global stiffness matrix has to be \textit{positive definite} in order for the discrete chemo-mechanical state $\Bd_1^*$ to be locally stable. As pointed out in the previous subsection, positive-definiteness of $\BK$ is a direct consequence of the space-time-discrete minimization principle (\ref{eq:2.24}). Therefore, the body becomes structurally unstable and bifurcates to an alternate buckled configuration at the instant when $\BK$ loses its positive-definiteness. This yields the eigenvalue problem for the matrix $\BK$ given by
\begin{equation}
\label{eq:2.39}
\left[\BK-\lambda\textbf{1}\right]\delta\Bd = \textbf{0}\,. 
\end{equation}  
Since all the eigenvalues of a real symmetric positive-definite matrix are positive, the bifurcation point occurs at the instant when the smallest eigenvalue $\lambda_{min}$ becomes negative 
\begin{equation}
\label{eq:2.40}
\boxed{
\lambda_{min}
\begin{cases}
> 0 \,,\quad \text{structurally stable chemo-mechanical state} \\[0.75ex]
\leq 0 \,,\quad \text{structurally unstable chemo-mechanical state}
\end{cases}}\,\,.
\end{equation}
The eigenvector $\delta\Bd$ corresponding to the critical eigenvalue $\lambda_{min}$ offers a representation of the buckled configuration of the body that minimizes its overall potential energy.

\section{Transient stability analysis of composite hydrogel structures}
In this section, we carry out the numerical implementation of the discrete minimization principle and the associated structural stability analysis to model the swelling-induced buckling of a pair of composite hydrogel structures under geometrical constraints. For numerically implementing the space-time-discrete minimization principle
\eqref{eq:2.24}
we employ a conforming finite-element design according to \citet{raviart1977}, \citet{brezzifortin} and \citet{teichtmeister2019}, see also Appendix~A.

\subsection{A constitutive model describing hydrogels}

The total strain-energy density $\hat{\psi}$ in the hydrogel includes three individual contributions and is given by
\begin{equation}
	\label{eq:3.1}
	\hat{\psi}(\BF,s)=\hat{\psi}_{mech}(\BF)+\hat{\psi}_{chem}(s)+\hat{\psi}_{coup}(J,s)\,.
\end{equation}     
The mechanical part $\hat{\psi}_{mech}(\BF)$, representing the elastic strain energy stored in the polymer chains of the gel, is modeled using a standard neo-Hookean energy function. For the chemical part $\hat{\psi}_{chem}(s)$, we adopt a Flory--Rehner-type function (refer \citet{flory1943}), which is based on the statistical thermodynamics of fluid-polymer interactions. It is assumed that in the dry state, the gel has negligible pore spaces and that the chemical bonds in the gel are strong enough to resist molecular deformations under external forces. As a result, the local volumetric strain in the gel is solely due to the diffusion of fluid molecules and this gives rise to the following molecular incompressibility constraint    
\begin{equation}
	\label{eq:3.2}
	J = 1+s\,,
\end{equation}
refer \citet{kang2010}, \citet{hong2009,hong2008}. The coupling term in (\ref{eq:3.1}) is a penalty enforcement of (\ref{eq:3.2}). For the dissipation potential $\hat{\phi}$, we consider a convex homogeneous function in $\mathbb{H}$ of degree two at a known chemo-mechanical state given by the right Cauchy--Green tensor $\BC_n$ and fluid concentration $s_n$ at time $t_n$. The logarithmic nature of the Flory--Rehner function for $\hat{\psi}_{chem}(s)$ results in a singularity at the dry state characterized by $s=0$. This causes problems in the numerical simulations when the dry state of the gel is considered as the reference state. To avoid this, we follow the procedure described in \citet{hong2009} and consider a \textit{preswollen} stress-free state as the reference configuration. To this end, the deformation gradient $\BF_d$, which maps the dry state of the gel to its current deformed configuration, is split multiplicatively as
\begin{equation}
	\label{eq:3.3}
	\BF_d=\BF\BF_0\,,
\end{equation}
where $\BF_0=J_0^{1/3}\textbf{1}$ represents an isotropic deformation of the dry gel resulting in a new stress-free preswollen reference state and $\BF$ is the actual deformation gradient relative to the preswollen state. The free-energy and dissipation potential functions are then suitably transformed using the Jacobian $J_0$ in order to comply with the new reference configuration of the gel. In the end, we have the following closed-form expressions for $\hat{\psi}$ and $\hat{\phi}$ given by   
\begin{equation}
	\label{eq:3.4}
	\begin{aligned}
	\hspace{-1mm}\hat{\psi}=\frac{\gamma}{2J_0}[\, J_0^{2/3}\BF:\BF-3-2\,\text{ln}(JJ_0)\,]&+
	\frac{\alpha}{J_0}[\,s\,\text{ln}(\frac{s}{1+s})+\frac{\chi s}{1+s}\, ] + \frac{\epsilon}{2J_0}(JJ_0-1-s)^2\\[1ex] \hspace{-2mm} \text{and} \qquad
	\hat{\phi}&=\frac{1}{2J_0^{1/3}M s_n}\BC_n:(\mathbb{H} \otimes \mathbb{H})\,,
	\end{aligned}
\end{equation}
respectively, refer \citet{boeger2017}. Using the stress-free condition of the preswollen reference state, the initial fluid concentration is determined to be
\begin{equation}
	\label{eq:3.5}
	s_0=\frac{\gamma}{\epsilon}(J_0^{-1/3}-\frac{1}{J_0})+J_0-1
\end{equation}
and using the consitutive equation (\ref{eq:2.6})\textsubscript{2}, the initial chemical potential of the gel is obtained as
\begin{equation}
	\label{eq:3.6}
	\mu_0=-\frac{\epsilon}{J_0}[J_0-1-s_0]+\frac{\alpha}{J_0}[\,\text{ln}(\frac{s_0}{1+s_0})+\frac{1}{1+s_0}+\frac{\chi}{(1+s_0)^2}\,]\,.
\end{equation}
The material parameters of the constitutive model (\ref{eq:3.4}) are described in Table \ref{tab:matpara}.  
\begin{table}[t]
	\captionsetup{width=0.9\textwidth}
	\caption{Material parameters describing the consititutive response of hydrogels}
	\footnotesize
	\centering
	\renewcommand{\arraystretch}{1.5}
	\begin{tabular*}{0.7\textwidth}{l@{\extracolsep{\fill}}lll}
		\hline
		No. & Parameter & Description & Units \\
		\hline 
		1. & $\gamma$ & Shear modulus & N/mm\textsuperscript{2} \\
		2. & $\alpha$ & Mixing modulus & N/mm\textsuperscript{2} \\
		3. & $\epsilon$ & Penalty parameter & N/mm\textsuperscript{2} \\
		4. & $M$ & Mobility parameter & mm\textsuperscript{4}/Ns\\
		5. & $J_0$ & Preswelling factor & - \\
		6. & $\chi$ & Interaction parameter & - \\
		\hline  
	\end{tabular*}
	\label{tab:matpara}
\end{table}  

\subsection{Surface wrinkling of flat hydrogel bilayers}

 \begin{figure}[h]
	\centering
	\scriptsize
	\psfrag{w}[c][c]{$w$}
	\psfrag{H}[c][c]{$H$}
	\psfrag{L}[c][c]{$L$}
	\psfrag{mucon}[l][l]{$\mu_0 \to \bar{\mu}$}
	\psfrag{aa}[r][r]{a)}
	\psfrag{bb}[r][r]{b)}
	\psfrag{t}[c][c]{$t/$s}
	\psfrag{zero}[c][c]{$0$}
	\psfrag{one}[c][c]{$1$}
	\psfrag{mu0}[r][r]{$\mu_0$}
	\psfrag{mubaar}[r][r]{$\bar{\mu}$}
	\psfrag{mu}[r][r]{$\mu$}
	\psfrag{x}[c][c]{$x$}
	\psfrag{y}[c][c]{$y$}
	\psfrag{film}[c][c]{\footnotesize film}
	\psfrag{sub}[c][c]{\footnotesize substrate}
	\includegraphics[width=0.9\textwidth]{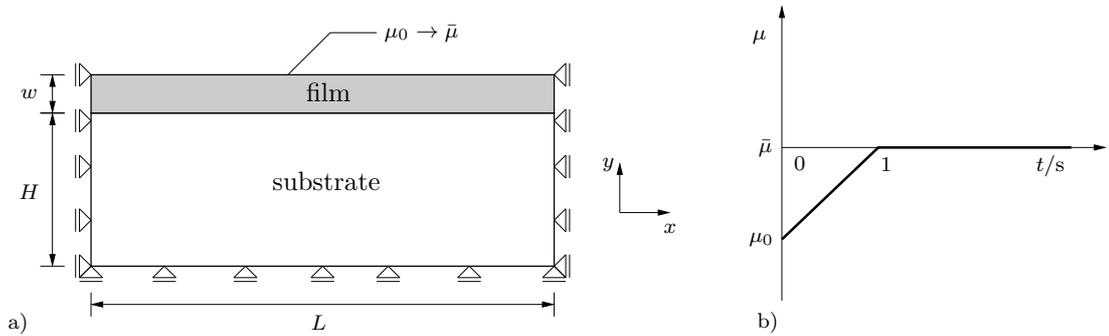}
	\captionsetup{width=0.95\textwidth}
	\caption{\footnotesize \textit{Geometric description of the hydrogel bilayer and the applied load profile.} a) Fluid diffuses into the bilayer from the top of the film and is allowed to accumulate inside by making the remaining edges of the bilayer impermeable. Suitable mechanical boundary conditions are applied so that the bilayer is allowed to swell only along the direction normal to the film surface. The dimensions of the substrate are taken to be $H=0.5\,$mm and $L=2\,$mm for all the analyses. b) Chemical potential on the film surface is increased linearly from its initial value $\mu_0$ to that of the diffusing fluid $\bar{\mu}=0$ in $1\,$s and then held constant. The resulting gradient in chemical potential along the bilayer thickness drives the diffusion process.}
	\label{fig:bi_geo}
\end{figure} 
In this subsection, we investigate the critical conditions for the onset of sinusoidal wrinkles on the surface of flat hydrogel bilayers that are subjected to geometrically constrained swelling. Fig.~\ref{fig:bi_geo}a depicts the geometry of the hydrogel bilayer under consideration along with the applied boundary conditions. It constitutes a stiff film of initial thickness $w$ perfectly bonded to a soft compliant substrate of initial thickness $H=0.5$ mm. Both layers are of equal length $L=2$ mm. Longitudinal expansion of the bilayer and the vertical displacement of the bottom edge are constrained. The specimen is loaded by increasing the chemical potential on the film surface from its initial value $\mu_0$ to that of the diffusing fluid $\bar{\mu}=0$ in $1\,$s, as shown in Fig.~\ref{fig:bi_geo}b. This establishes a gradient in chemical potential along the thickness, causing fluid to  diffuse through the film. The remaining edges are made impermeable to fluid outflux so that the diffusing fluid is allowed to accumulate inside the bulk of the bilayer. As a result of fluid accumulation, the specimen begins to swell but only along the outward normal to the film surface due to the prescribed mechanical boundary conditions. The longitudinal constraint generates a compressive stress in the structure, which increases progressively with increasing amount of fluid accumulation. Once the compressive stress reaches a critical value, sinusoidal wrinkles are observed on the film surface which are seen as manifestations of a structural instability activated in the system to relieve some of the excess compression.
\begin{table}[h]
	\captionsetup{width=0.8\textwidth}
	\caption{\footnotesize{The set of material parameters that are kept constant for all the analyses on the bilayer system. The values are identical for both the film and the substrate.}}
	\footnotesize
	\centering
	\renewcommand{\arraystretch}{1.6}
	\begin{tabular*}{0.7\textwidth}{l@{\extracolsep{\fill}}lll}
		\hline
		Parameter & Description & Value & Units \\ 
		\hline   
		$\alpha$ & Mixing modulus & $24.2$ & MPa\\
		$\epsilon$ & Penalty parameter & $10$ & MPa\\
		$M$ & Mobility parameter & $10^{-4}$ & mm\textsuperscript{4}/Ns\\
		\hline
	\end{tabular*}
	
	\label{tab:bilayerpara}
\end{table} 
 
We characterize the structural instability using the critical growth $g_c$ and the number of wrinkles $N_c$. The former is defined as the increase in total thickness of the bilayer for which wrinkles start to emerge on the film surface. In what follows, we study the dependence of $g_c$ and $N_c$ on the geometrical and material parameters of the bilayer\footnotemark[1]. Table \ref{tab:bilayerpara} summarizes the values of those material parameters that are held constant for all the forthcoming investigations and these values are identical for both the film and the substrate, refer \citet{boeger2017}. Henceforth, we adopt the superscript $f$ for the quantities of the film and the superscript $s$ for those of the substrate.  
\footnotetext[1]{At this point, it is worth mentioning that the critical buckling loads estimated in terms of the chemical potential at the film surface are not accurate quantities. In the minimization-based formulation, the chemical-potential field appears as a Neumann variable (refer Eq.~(\ref{eq:2.10})) and therefore, any boundary condition imposed on it is only satisfied in a \textit{weak sense}. Furthermore, the chemical potentials are computed as part of the driving-force array $\mathfrak{D}$ at the Gauss-quadrature points of Q\textsubscript{1}RT\textsubscript{0} elements and these are then \textit{projected} onto the nodal points during the post-processing stage. A combination of these weak-satisfaction and projection errors render the chemical-potential values unreliable for estimating buckling loads. As a result, we have chosen the specimen growth as the critical parameter to characterize the buckling phenomenon. In addition to being more accurate, it is also a quantity that can be easily monitored in an experimental setting.}

\begin{figure}[t]
	\centering
	\footnotesize
	\psfrag{g}[c][c]{$g_c\,/\mu$m}
	\psfrag{w}[c][c]{$w/$mm}
	\psfrag{aa}[c][c]{a)}
	\psfrag{bb}[c][c]{b)}
	\psfrag{N}[c][c]{$N_c$}
	\psfrag{1}[r][r]{$N=1$}
	\psfrag{1.5}[r][r]{$N=1.5$}
	\psfrag{2}[r][r]{$N=2$}
	\includegraphics[width=0.95\textwidth]{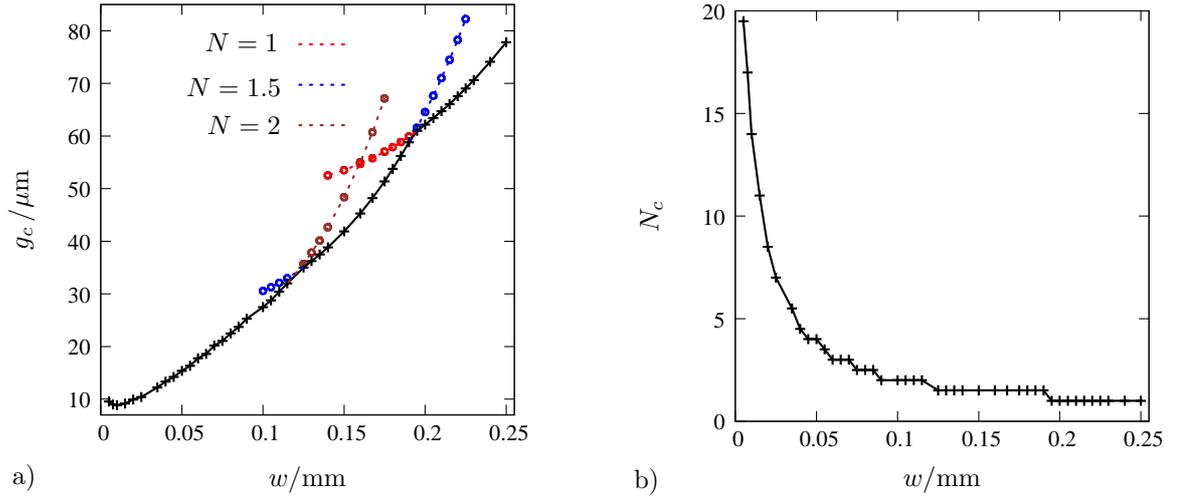}
	\captionsetup{width=0.95\textwidth}
	\caption{\footnotesize \textit{Influence of film thickness on the critical buckling characteristics.} a) Critical growth $g_c$ and b) wrinkle count $N_c$ as functions of the film thickness $w$ for a bilayer having modulus ratio $\gamma^f/\gamma^s=8$, preswelling factor $J_0^{f,s}=1.01$ and interaction parameter $\chi^{f,s}=0.1$. The colored curves $N=\left\lbrace1,1.5,2 \right\rbrace $ represent the extensions of the respective isomodes. Intersection of two successive isomodes indicates a shift in buckling mode due to the existence of an alternate buckled state at a lower critical load. }
	\label{fig:bi_film}
\end{figure}

\subsubsection{Influence of film thickness on the buckling characteristics}

Fig.~\ref{fig:bi_film} illustrates the influence of film thickness $w$ on the critical buckling characteristics of a hydrogel bilayer having a shear-modulus ratio $\gamma^f/\gamma^s=8$. For the preswelling factor and interaction parameter, the values $J_0^{f,s}=1.01$ and $\chi^{f,s}=0.1$ are used. The critical growth $g_c$ is found to increase monotonically with the film thickness. This is to be expected as a bilayer having a thicker film requires a greater amount of accumulated fluid to generate enough longitudinal compression that can trigger the bifurcation mode. A similar monotonic trend has been illustrated by \citet{jin2015} for the uniaxial homogeneous compression of a soft neo-Hookean bilayer, where the critical longitudinal compressive strain is plotted against the substrate-film thickness ratio. Furthermore, the curve showing $g_c$ as a function of $w$ in Fig.~\ref{fig:bi_film}a is the envelope of a family of \textit{isomodes.} An isomode is a curve showing the variation of $g_c$ with respect to $w$ for a specific mode of buckling with wrinkle count $N$. The colored curves in Fig.~\ref{fig:bi_film}a represent the extensions of isomodes corresponding to wrinkle counts mentioned in the plot. An intersection of two successive isomodes indicates a shift in the buckling mode and the primary reason for this shift is the existence of an alternate buckled state at a \textit{lower} critical load. For instance, if the mode $N=1$ had not existed, the curve would have continued along the extension of $N=1.5$ shown in blue. Such shifts between buckling modes happen all along the $g_c$ curve, however, the intersection of isomodes has not been indicated explicitly at lower film thicknesses as they are much more closely spaced. This behaviour bears similarities to the plot showing the dependence of buckling load of a plate on its aspect ratio during uniaxial compression, wherein a sudden decrease in the buckling load is observed due to the shift in mode number, refer \citet{reddy2006}. Increasing the film thickness also leads to fewer number of wrinkles appearing in the buckled state, as seen in Fig.~\ref{fig:bi_film}b. This is attributed to a decrease in film flexibility at higher thicknesses making it more difficult to be bent into wrinkles.
 
\subsubsection{Influence of modulus ratio on the buckling characteristics}
 
\begin{figure}[t]
        \footnotesize
	\centering
	\psfrag{g}[c][c]{$g_c\,/\mu$m}
	\psfrag{MR}[c][c]{$\gamma^f/\gamma^s$}
	\psfrag{aa}[c][c]{a)}
	\psfrag{bb}[c][c]{b)}
	\psfrag{WN}[c][c]{$N_c\,$}
	\includegraphics[width=0.95\textwidth]{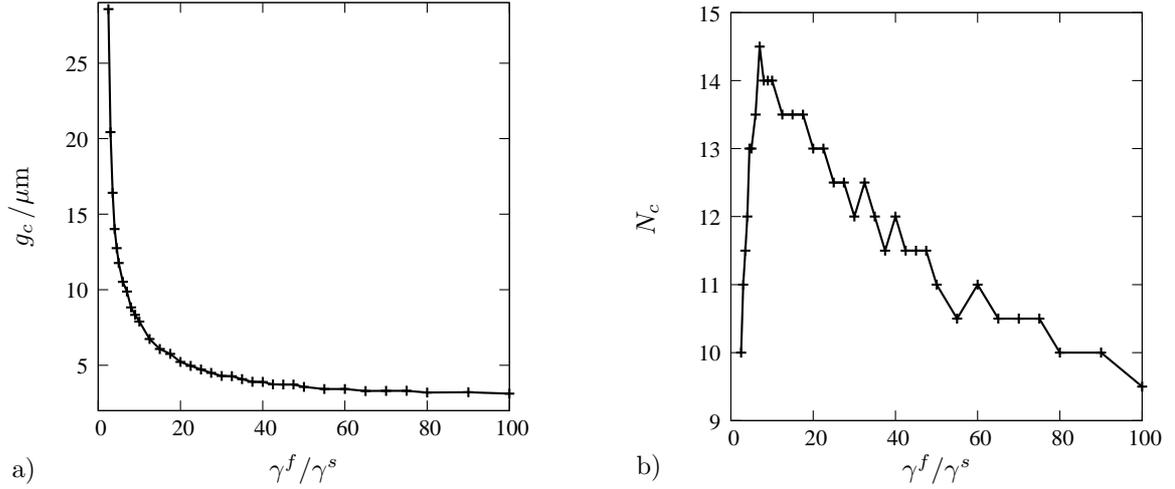}
	\captionsetup{width=0.95\textwidth}
	\caption{\footnotesize \textit{Effect of shear-modulus ratio on the critical buckling characteristics.} Variations in a) critical growth $g_c$ and b) wrinkle count $N_c$ with respect to the shear-modulus ratio $\gamma^f/\gamma^s$ of a hydrogel bilayer with film thickness $w=0.01\,$mm. The material parameters $J_0^{f,s}=1.01$ and $\chi^{f,s}=0.1$ are used. $g_c$ is observed to be a decreasing function of $\gamma^f/\gamma^s$, whereas $N_c$ shows an increasing-decreasing trend.} 	
	\label{fig:bi-modrat}
\end{figure}

Fig.~\ref{fig:bi-modrat} describes the critical buckling characteristics as functions of the shear modulus ratio for a bilayer having film thickness $w=0.01\,$mm, preswelling factor $J_0^{f,s}=1.01$ and interaction parameter $\chi^{f,s}=0.1$. The critical growth is found to decrease as the film gets stiffer. 
\begin{figure}[b]
	\centering
	\psfrag{w}[c][c]{$w/$mm}
	\psfrag{MR}[c][c]{$\gamma^f/\gamma^s\,$}
	\includegraphics[width=\textwidth]{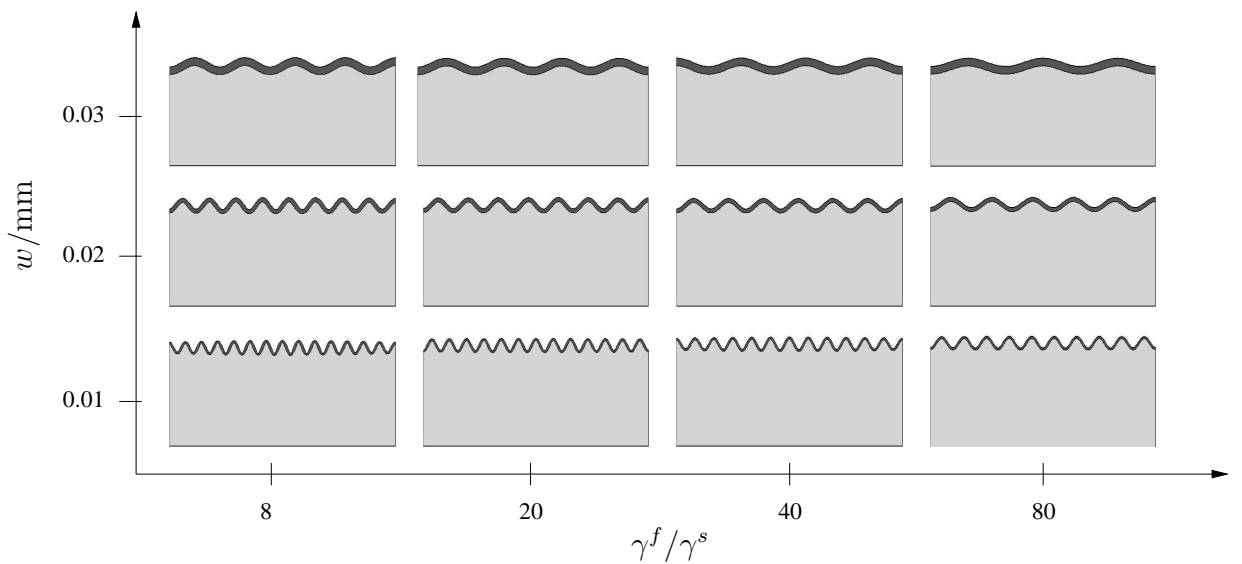}
	\captionsetup{width=\textwidth}
	\caption{\footnotesize \textit{Critical buckling modes of the bilayer for selected film thicknesses and shear-modulus ratios.} For a given modulus ratio, fewer wrinkles are seen with increasing film thickness as described in Fig.~\ref{fig:bi_film}b. Furthermore, the number of wrinkles can be seen to decrease with increasing modulus ratio for a fixed bilayer geometry as illustrated in Fig.~\ref{fig:bi-modrat}b.}
	\label{fig:bi_boxplot}
\end{figure}
The reason for this behaviour is that for a given volume of accumulated fluid, a stiffer film experiences a greater longitudinal compressive stress, leading to the critical state of compression being reached for lower values of growth. For moderate modulus ratios ($\gamma^f/\gamma^s<10$), there is a sharp decrease in $g_c$ with increasing $\gamma^f/\gamma^s$ whereas the decrease is only marginal at higher modulus ratios. In the latter regime, the number of wrinkles in the critical buckling mode is found to decrease with increasing shear-modulus ratio as seen in Fig.~\ref{fig:bi-modrat}b. This is due to the fact that a film with a higher shear modulus is more resistant to being bent into wrinkles. However, when the two layers have comparable moduli ($\gamma^f/\gamma^s<10$), there is an increase in wrinkle count with respect to the modulus ratio. Such an increasing-decreasing trend in the mode number has also been reported by \citet{cao2012} for the buckling of neo-Hookean bilayers under plane-strain compression. Interestingly, the interval of comparable film and substrate moduli for which $N_c$ increases with $\gamma^f/\gamma^s$ is quite identical to the range of modulus ratios for which $g_c$ shows a sharp decrease as seen in Fig.~\ref{fig:bi-modrat}a. The critical buckling modes of the bilayer for some selected film thicknesses and shear-modulus ratios are shown in Fig.~\ref{fig:bi_boxplot}. The trends observed in Fig.~\ref{fig:bi_film}b and Fig.~\ref{fig:bi-modrat}b are clearly evident.
\begin{figure}[t]
	\centering
	\scriptsize
	\psfrag{g}[c][c]{$g_c\,/\mu$m}
	\psfrag{0.1}[r][r]{ $\chi^{f,s}=0.1$}
	\psfrag{0.8}[r][r]{ $\chi^{f,s}=0.8$}
	\psfrag{J0}[c][c]{$J_0^{f,s}\,\,$}
	\psfrag{aa}[r][r]{a)}
	\psfrag{bb}[r][r]{b)}
	\psfrag{cc}[r][r]{c)}
	\psfrag{N}[c][c]{$N_c\,\,$}
	\psfrag{all}[r][r]{ $\chi^{f,s}=\left\lbrace0.1,0.8 \right\rbrace $}
	\psfrag{mu}[c][c]{$\mu_c\,\,$[MPa]}
	\psfrag{initial0.1}[l][l]{ $\mu_0^{\left(\chi=0.1 \right)} $}
	\psfrag{initial0.8}[l][l]{ $\mu_0^{\left(\chi=0.8 \right)} $}
	\includegraphics[width=0.95\textwidth]{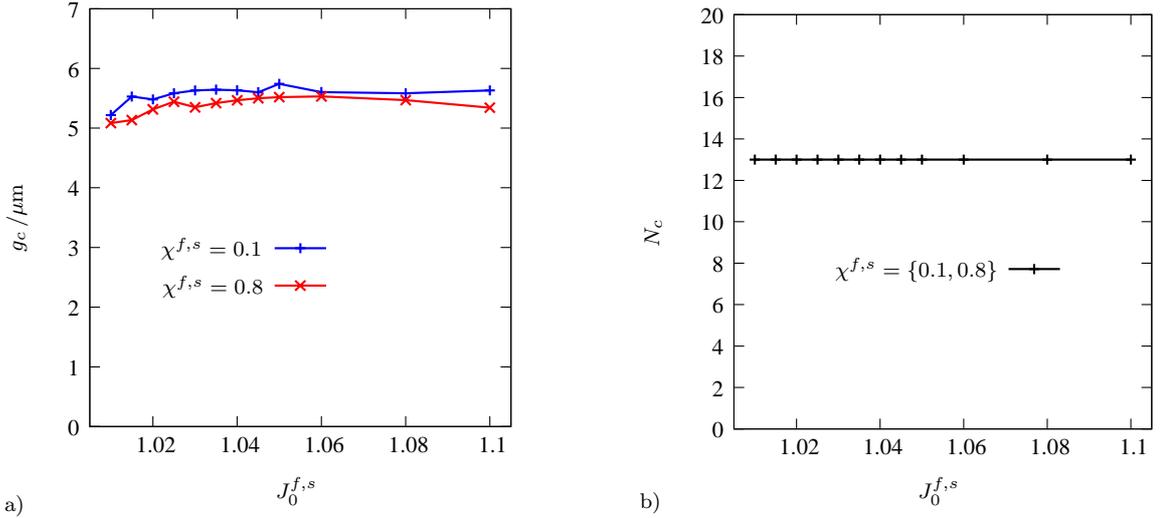}
	\captionsetup{width=0.95\textwidth}
	\caption{\footnotesize \textit{Critical buckling characterisitics of the bilayer as functions of the preswelling factor and interaction parameter.} Dependence of a) critical growth $g_c$ and b) wrinkle count $N_c$ on the initial preswelling factor $J_0^{f,s}$ for interaction parameters $\chi^{f,s}=\left\lbrace0.1,0.8 \right\rbrace $. A bilayer with film thickness $w=0.01\,$mm and modulus ratio $\gamma^f/\gamma^s=20$ is chosen for this study. The buckling characteristics $g_c$ and $N_c$ are found to be largely unaffected by both $J_0^{f,s}$ and $\chi^{f,s}$.}
	\label{fig:bi_preswelling}
\end{figure}

\subsubsection{Influence of preswelling factor and interaction parameter on the buckling characteristics}

The dependence of critical buckling characteristics on the preswelling factor and interaction parameter is illustrated in Fig.~\ref{fig:bi_preswelling} for a bilayer having film thickness $w=0.01\,$mm and shear-modulus ratio $\gamma^f/\gamma^s=20$. Compared to the dependences on film thickness and modulus ratio, the critical growth remains largely unchanged with respect to both $J_0$ and $\chi$, the variations being within $1\,$$\mu$m. From Fig.~\ref{fig:bi_preswelling}b, we infer that $N_c$ is also independent of the initial preswelling factor and interaction parameter. This behaviour is to be expected as the transformation to a preswollen initial configuration does not result in any prestress in the bilayer. As a result, bifurcation happens at the same critical growth irrespective of the extent of initial preswelling, leading to a buckled configuration having the same number of wrinkles since neither the film thickness nor the modulus ratio is altered. 

\subsection{Surface wrinkling of bilayered hydrogel tubes}

\begin{figure}[t]
	\centering
	\footnotesize
	\psfrag{aa}[c][c]{a)}
	\psfrag{bb}[c][c]{b)}
	\psfrag{H}[c][c]{\scriptsize $H$}
	\psfrag{w}[c][c]{\scriptsize $w$}
	\psfrag{r}[c][c]{\scriptsize $r$}
	\psfrag{x}[c][c]{$x$}
	\psfrag{y}[c][c]{$y$}
	\psfrag{mucon}[l][l]{\scriptsize $\mu_0 \to \bar{\mu} $}
	\psfrag{film}[l][l]{\scriptsize film}
	\psfrag{sub}[l][l]{\scriptsize substrate}
	\includegraphics[width=0.9\textwidth]{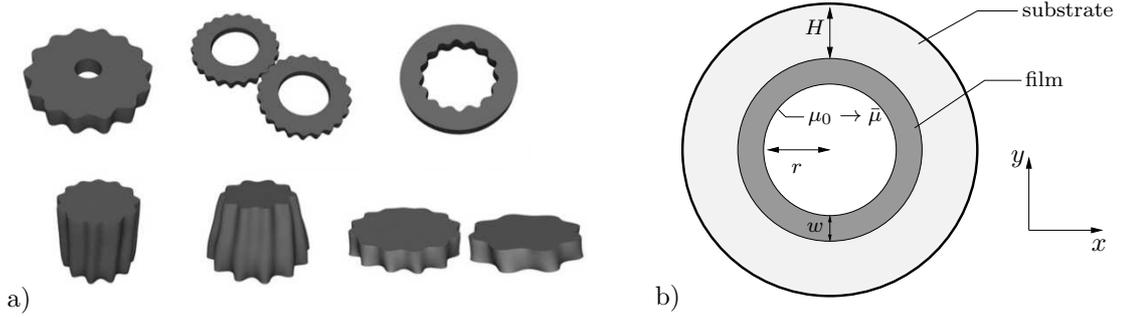}
	\captionsetup{width=0.95\textwidth}
	\caption{\footnotesize \textit{Cross-sectional geometry of the hydrogel tube together with the boundary conditions.} a) Gear profiles generated through the controlled buckling of cylindrical film-substrate bilayers as described in \citet{yin2009} and b) cross-sectional geometry of the representative hydrogel tube. A stiff film of thickness $w$ is deposited on the inner wall of a soft compliant cylindrical substrate of thickness $H$. The outer surface of the substrate, having a fixed radius of $1\,$mm, has a rigid and impermeable wall (represented by a bold outer circle) that prevents radial expansion. Fluid is allowed to diffuse into the film by increasing the surface chemical potential from its initial value $\mu_0$ to that of the fluid $\bar{\mu}=0$ in $1\,$s.  }
	\label{fig:tu_geo}
\end{figure}

In this subsection, we study the swelling-induced sinusoidal wrinkling on the inner surface of representative bilayered hydrogel tubes having rigid outer walls. \citet{yin2009} have demonstrated the fabrication of soft microgears by a mechanical self-assembly process wherein complex gear profiles, with desired number of teeth and amplitude, are generated purely by exploiting the buckling modes of cylindrical film-substrate systems, see Fig.~\ref{fig:tu_geo}a. Inner-surface wrinkling of constrained soft materials is also observed as a common natural phenomenon. Mucosal tissues, that form the inner lining of tubular organs such as the oesophagus, are subjected to growth under the constraint of stiff muscular layers. The resulting compressive residual stresses destabilize the tissues, leading to the formation of wrinkles on the inner surface. This has been investigated in detail by \citet{li2011}. For our analyses, we assume plane-strain conditions and focus on a specific cross-section of the cylindrical tube. The geometry and boundary conditions are described in Fig.~\ref{fig:tu_geo}b.

The tube has an outer radius of $1\,$mm. A stiff film of thickness $w$ is bound to the inner surface of a cylindrical compliant substrate of thickness $H$. The initial chemical potential $\mu_0$ on the film is increased to a value $\bar{\mu}=0$ in $1\,$s similar to the profile shown in Fig.~\ref{fig:bi_geo}b, allowing fluid to diffuse radially into the bilayer. A rigid impermeable wall surrounding the outer surface of the substrate prevents radial expansion and fluid outflux. Under these constraints, the film-substrate system begins to swell radially inwards, shrinking the radius $r$ with the passage of time. As a result, the circumferential compressive stress in the film increases continuously and once a certain critical state of compression is attained, a bifurcation mode is activated resulting in the formation of sinusoidal wrinkles on the film surface. In the following, we investigate the dependence of the critical buckling characteristics, namely the critical radial growth $g_c$ and the number of wrinkles $N_c$ along the inner circumference, on the geometry and material parameters of the bilayered tube. In addition to the values of material parameters summarized in Table \ref{tab:bilayerpara}, we also set fixed values for the preswelling factor $J_0^{f,s}=1.01$ and interaction parameter $\chi^{f,s}=0.1$ for all subsequent analyses.

\subsubsection{Dependence of buckling characteristics on the film thickness}

\begin{figure}[t]
	\centering
	\footnotesize
	\psfrag{g}[c][c]{$g_c\,/\mu$m}
	\psfrag{w}[c][c]{$w\,/$mm}
	\psfrag{aa}[c][c]{a)}
	\psfrag{bb}[c][c]{b)}
	\psfrag{N}[c][c]{$N_c$}
	\includegraphics[width=0.95\textwidth]{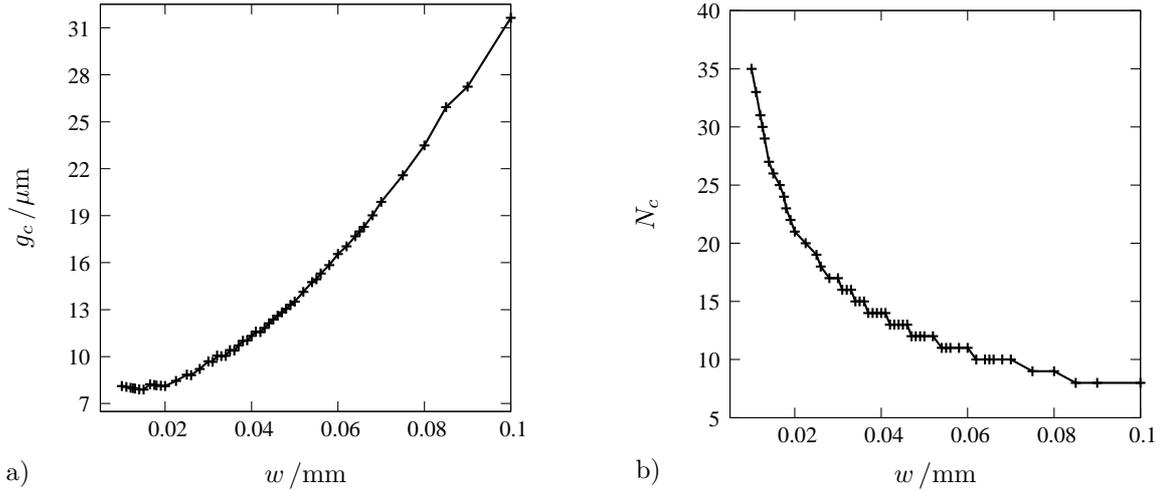}
	\captionsetup{width=0.95\textwidth}
	\caption{\footnotesize \textit{Effect of film thickness on the critical buckling characteristics of the hydrogel tube.} a) Critical growth $g_c$ and b) wrinkle count $N_c$ as functions of the film thickness $w$ for a tube having a substrate thickness $H=0.2\,$mm and shear modulus ratio $\gamma^f/\gamma^s=10.$ We find the trends to be similar to the ones seen in Fig.~\ref{fig:bi_film}a) for flat bilayers. However, in the present case of hydrogel tubes, the isomodes are much closer to each other and as a result, shifts between the buckling modes are not explicitly seen.}
	\label{fig:tu_film}
\end{figure}

Fig.~\ref{fig:tu_film} illustrates the effect of varying the film thickness $w$ on the critical buckling characteristics of a hydrogel tube, having substrate thickness $H=0.2\,$mm and shear-modulus ratio $\gamma^f/\gamma^s=10$. The critical growth $g_c$ is found to increase monotonically with the film thickness of the bilayered tube. This can be explained by the fact that a thicker film requires a larger circumferential compressive stress to destabilize the system and to achieve this, a greater amount of swelling-induced growth is needed. Similar to Fig.~\ref{fig:bi_film}a, the curve in Fig.\ref{fig:tu_film}a is the envelope of a family of isomodes. In this case, shifts between the successive isomodes are not seen explicitly as they are a lot closer to each other. In other words, the surface morphology of the buckled configuration is more sensitive to the film thickness than in the case of a flat hydrogel bilayer. As seen in Fig.~\ref{fig:tu_film}b, the number of wrinkles is found to decrease with increasing film thickness. This behaviour is attributed to the reduced flexibility of thicker films compared with thinner films, making it difficult for the former to be bent into wrinkles. This trend has also been observed in the works of \citet{li2011} and \citet{xie2014}, where the onset of wrinkling on the inner surface of cylindrical film-substrate systems described by an incompressible hyperelastic material law has been investigated both analytically and computationally.
 
\subsubsection{Dependence of buckling characteristics on the substrate thickness}

\begin{figure}[t]
	\centering
	\footnotesize
	\psfrag{g}[c][c]{$g_c\,/\mu$m}
	\psfrag{H}[c][c]{$H/$mm}
	\psfrag{N}[c][c]{$N_c$}
	\psfrag{aa}[c][c]{a)}
	\psfrag{bb}[c][c]{b)}
	\includegraphics[width=0.95\textwidth]{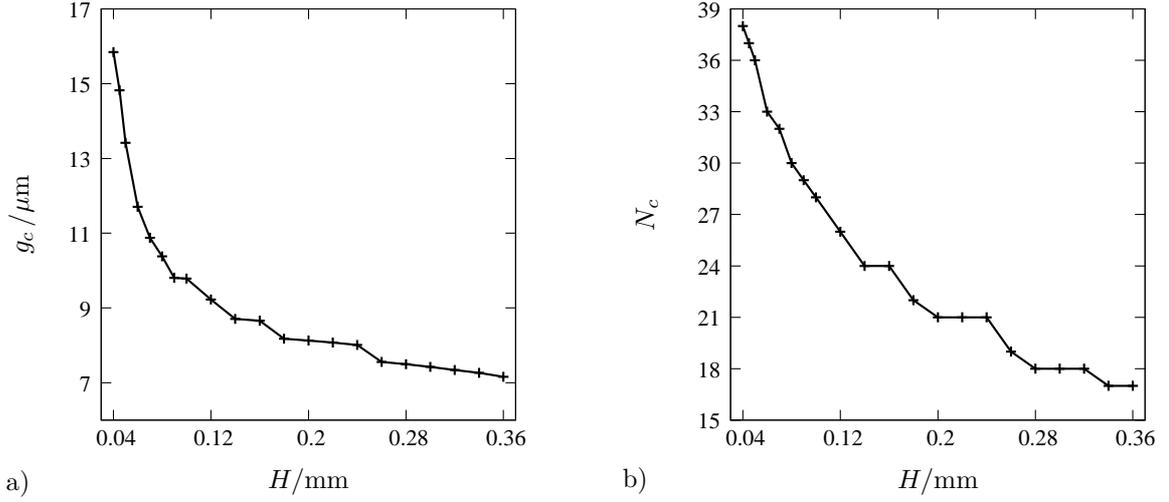}
	\captionsetup{width=0.95\textwidth}
	\caption{\footnotesize \textit{Effect of substrate thickness on the critical buckling characteristics of the tube.} a) Critical growth $g_c$ and b) wrinkle count $N_c$ as functions of the substrate thickness $H$ of a bilayered hydrogel tube having modulus ratio $\gamma^f/\gamma^s=10$ and film thickness $w=0.02\,$mm. Both $g_c$ and $N_c$ are found to decrease with increasing $H$.  }
	\label{fig:tu_substrate}
\end{figure}

The effects of changing the initial substrate thickness $H$ on the critical buckling characteristics of the tube are illustrated in Fig.~\ref{fig:tu_substrate}. For this study, a tube having $w=0.02\,$mm and $\gamma^f/\gamma^s=10$ is chosen. From Fig.~\ref{fig:tu_substrate}a, we observe that increasing the initial thickness of the substrate has a destabilizing effect on the tubular structure, leading to wrinkles being activated at lower values of critical growth. A similar trend has been observed in the work of \citet{moulton2011}, wherein the onset of growth-induced wrinkles in cylindrical hyperelastic tubes has been studied for various thicknesses by adopting an analytical approach. The onset of wrinkles can be considered as that instant when the substrate, under constrained growth, begins to experience shear deformation. The effective stiffness of the substrate against shear decreases with increasing initial thickness $H$. Therefore, for a given initial film thickness and loading rate, a thicker substrate tends to shear sooner. This explains why tubes with thicker substrates have lower values of critical growth. As seen in Fig.~\ref{fig:tu_substrate}b, tubes with thicker substrates tend to have fewer wrinkles in the critical buckling mode. As pointed out by \citet{li2011} in their study of wrinkling of soft tissues lining the inner surfaces of tubular biological organs, this is consistent with relevant practical observations wherein fewer wrinkles are seen on the surface of a stiff mucosal tissue when the softer submucosal layer surrounding it is thicker than normal.    

\subsubsection{Dependence of buckling characteristics on the modulus ratio}

\begin{figure}[b]
	\centering
	\footnotesize
	\psfrag{MR}[c][c]{$\gamma^f/\gamma^s\,$}
	\psfrag{g}[c][c]{$g_c\,/\mu$m}
	\psfrag{N}[c][c]{$N_c$}
	\psfrag{aa}[c][c]{a)}
	\psfrag{bb}[c][c]{b)}
	\includegraphics[width=0.95\textwidth]{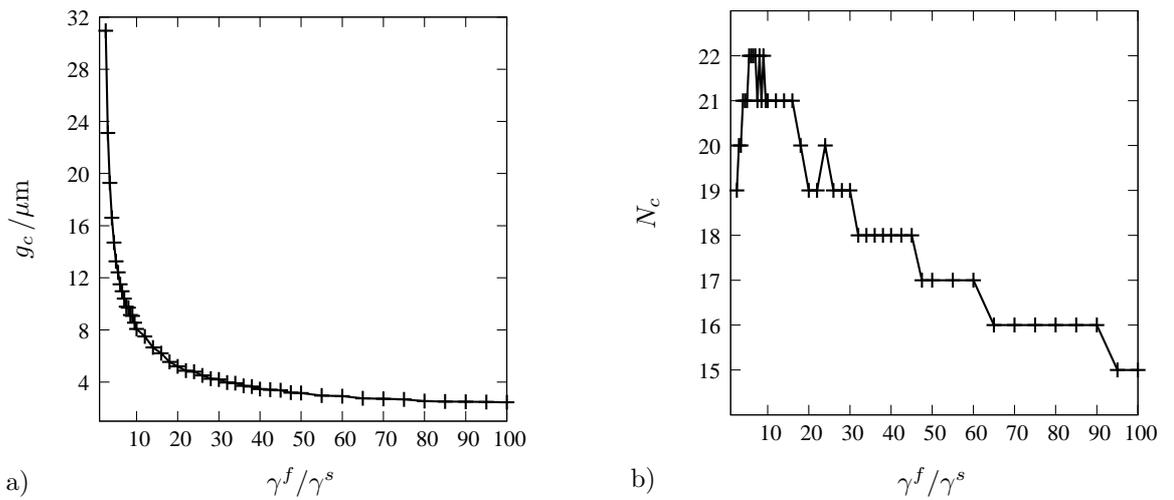}
	\captionsetup{width=0.95\textwidth}
	\caption{\footnotesize \textit{Effect of shear-modulus ratio on the critical buckling characteristics of the tube.} a) Critical growth $g_c$ and b) wrinkle count $N_c$ as functions of the shear-modulus ratio $\gamma^f/\gamma^s$ for a tube with film and substrate thicknesses $w=0.02\,$mm and $H=0.2\,$mm respectively. $g_c$ decreases sharply with increasing $\gamma^f/\gamma^s$ for lower modulus ratios. With $N_c$, we observe an increasing-decreasing trend similar to the case of a flat hydrogel bilayer.}
	\label{fig:tu_modrat}
\end{figure}	

Fig.~\ref{fig:tu_modrat} shows the dependence of the critical buckling characteristics on the shear-modulus ratio for a tube having film and substrate thicknesses $w=0.02\,$mm and $H=0.2\,$mm, respectively. 
\begin{figure}[t]
	\centering
	\psfrag{w}[c][c]{$w\,/$mm}
	\psfrag{MR}[c][c]{$\gamma^f/\gamma^s$}
	\includegraphics[width=0.8\textwidth]{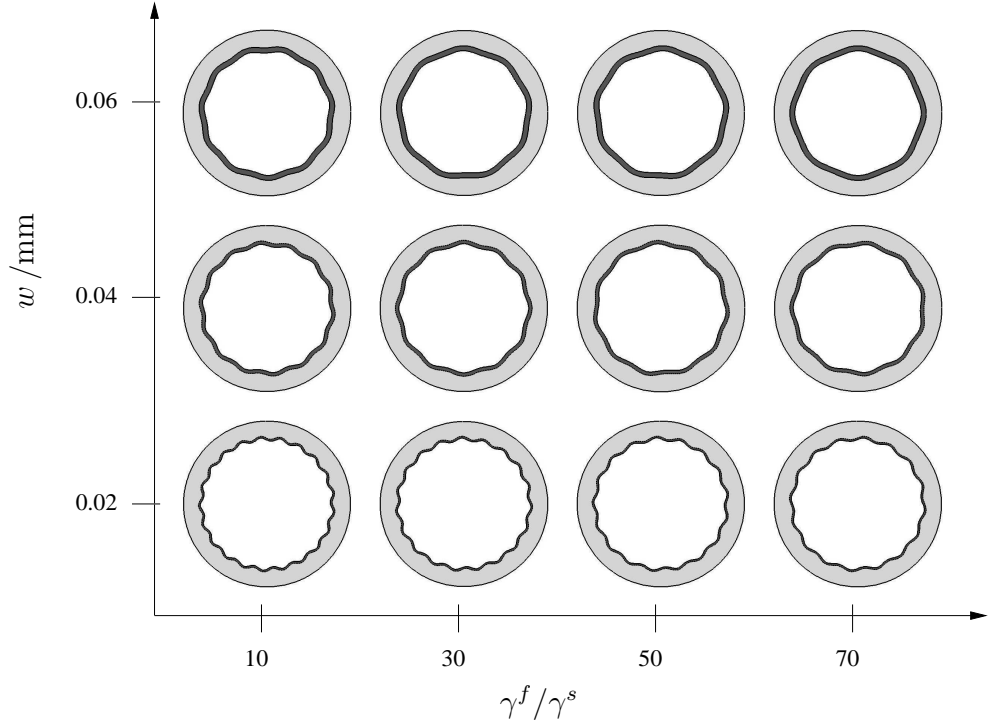}
	\captionsetup{width=0.9\textwidth}
	\caption{\footnotesize \textit{Critical buckling modes for some selected film thicknesses and modulus ratios.} It can be seen that increasing both $w$ and $\gamma^f/\gamma^s$ (beyond the range of comparable film and substrate moduli) leads to fewer wrinkles in the critical buckling mode.}
	\label{fig:tu_boxplot}
\end{figure}
\begin{figure}[b]
	\centering
	\footnotesize
	\psfrag{r}[c][c]{$r$}
	\psfrag{a}[c][c]{$b$}
	\psfrag{b}[c][c]{$a$}
	\psfrag{aa}[c][c]{a)}
	\psfrag{bb}[c][c]{b)}
	\psfrag{x}[c][c]{$x$}
	\psfrag{y}[c][c]{$y$}
	\includegraphics[width=0.9\textwidth]{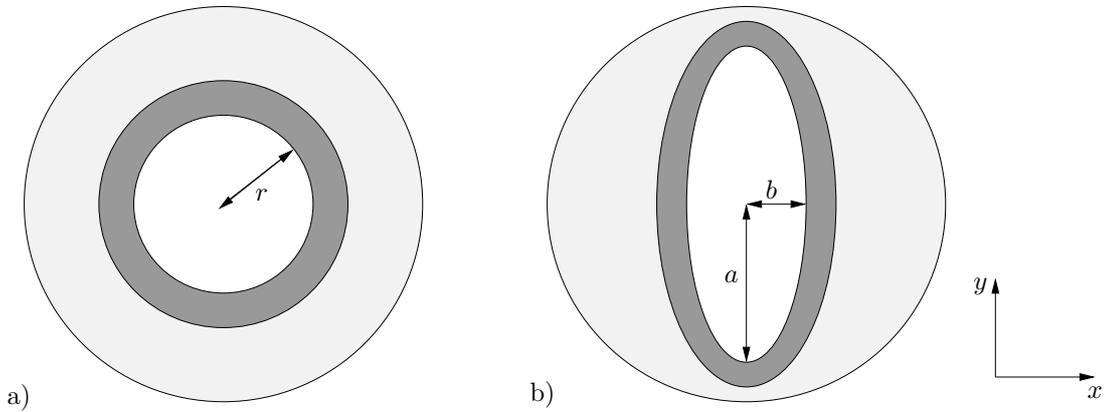}
	\captionsetup{width=0.9\textwidth}
	\caption{\footnotesize \textit{Perturbed cross-sectional geometry of the bilayered hydrogel tube.} a) The reference geometry of the cross-section having radius $r=0.78\,$mm and film thickness $w=0.02\,$mm. b) Tube with an elliptical passage having a semi-major axis $a$ and semi-minor axis $b$. The film thickness $w=0.02\,$mm is kept uniform along the elliptical circumference. For every perturbed cross-section, the relation $ab=r^2$ holds, i.e. the passage area is kept constant.}
	\label{fig:tu_pbgeo}
\end{figure}The critical growth is found to decrease monotonically with increasing modulus ratio. This is to be expected because for a given amount of radial growth in the specimen, a stiffer film experiences greater compressive stress and therefore, the critical state of compression is reached sooner at lower values of growth. In Fig.~\ref{fig:tu_modrat}a, we observe an initial sharp decrease in $g_c$ at moderate stiffness ratios ($\gamma^f/\gamma^s<20$) after which the decrease is much less prominent. Similar to the case of a flat hydrogel bilayer, there is an initial increase in $N_c$ with increasing $\gamma^f/\gamma^s$ when the film and substrate layers have comparable shear moduli, following which there is a subsequent decrease at higher modulus ratios. Similar trends have also been reported by \citet{jin2018} and \citet{li2020} in their investigations on wrinkling instability of cylindrical hyperelastic film-substrate systems. In the work of the latter, this rather abnormal trend in the mode number with respect to $\gamma^f/\gamma^s$ has been attributed to a transition from global buckling of the tubular structure to localized wrinkling on the film surface, i.e.\ when the film and the substrate have comparable shear moduli, the two layers deform in a more synchronized manner resulting in the tube buckling globally at the critical state of growth whereas at higher modulus ratios, the film reaches the critical state well ahead of the substrate. 
\begin{figure}[b]
	\centering
	\psfrag{1}[c][c]{$1$}
	\psfrag{1.1}[c][c]{$1.1$}
	\psfrag{1.2}[c][c]{$1.2$}
	\psfrag{1.3}[c][c]{$1.3$}
	\psfrag{1.4}[c][c]{$1.4$}
	\psfrag{1.5}[c][c]{$1.5$}
	\psfrag{2}[c][c]{$2$}
	\psfrag{4}[c][c]{$4$}
	\psfrag{6}[c][c]{$6$}
	\psfrag{8}[c][c]{$8$}
	\psfrag{10}[c][c]{$10$}
	\psfrag{12}[c][c]{$12$}
	\psfrag{14}[c][c]{$14$}
	\psfrag{gx}[c][c]{$g_{c,x}\,/\mu$m}
	\psfrag{a/b}[c][c]{$a/b$}
	\includegraphics[width=0.75\textwidth]{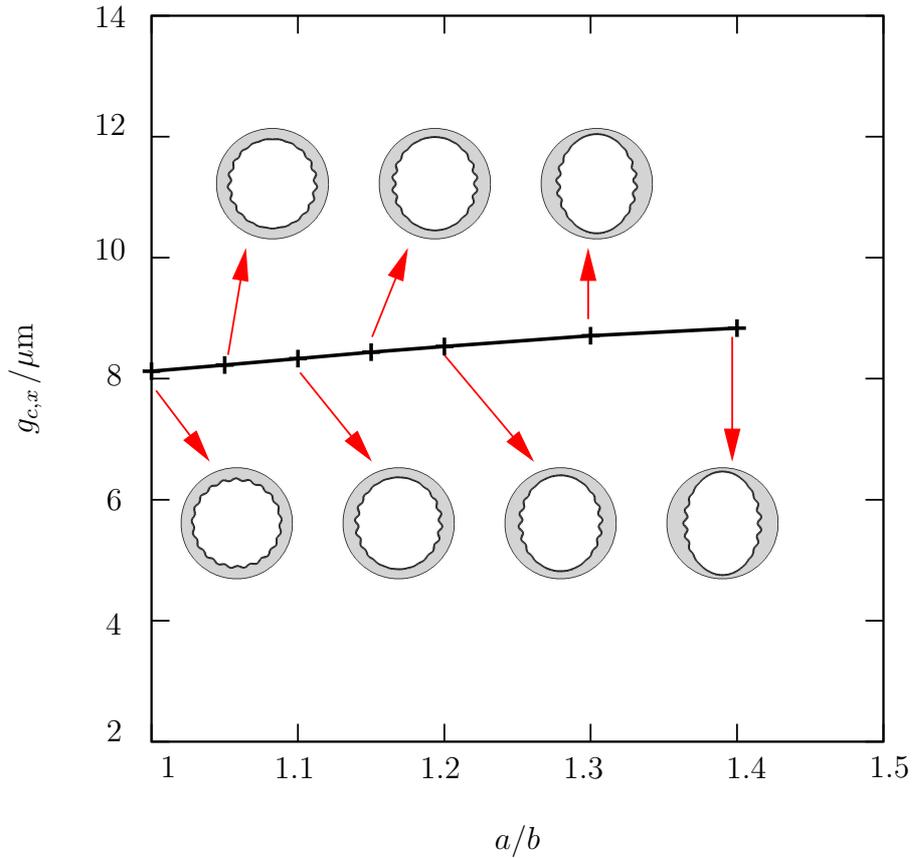}
	\captionsetup{width=0.95\textwidth}
	\caption{\footnotesize \textit{Influence of the passage shape on the critical growth and wrinkling pattern.} Critical growth $g_{c,x}$ along the minor axis as a function of the aspect ratio $a/b$ of the elliptical passage. For all perturbed geometries, the values $w=0.02\,$mm and $\gamma^f/\gamma^s=10$ are used. Although only a marginal increase in $g_{c,x}$ is observed , there is a significant change in the wrinkling pattern, with fewer or no wrinkles being formed in the regions where the substrate is thinner.}
	\label{fig:tu_psg}
\end{figure}

\subsubsection{Dependence of buckling characteristics on the passage shape} 

We now investigate the effects of perturbing the shape of the circular passage on the critical buckling characteristics of the tube. Fig.~\ref{fig:tu_pbgeo} shows the original geometry of the tube together with the perturbed cross-sectional geometry comprising an elliptical passage. The film thickness $w$ is kept uniform along the circumference of the ellipse, resulting in a non-uniform substrate thickness as seen in Fig.~\ref{fig:tu_pbgeo}b. The aspect ratio $a/b$ for each geometry is chosen such that the passage area always remains the same and matches with that of the reference tube in Fig.~\ref{fig:tu_pbgeo}a, i.e.\ the relation $ab=r^2$ holds for all perturbed cross-sections. The same boundary conditions and load profile are used as for the previous analyses. 

Owing to the elliptical shape of the passage, the radial growth is no longer axisymmetric, i.e. the magnitude of growth is observed to be greater along the minor axis than along the major axis. As a result, we focus on the critical growth along the minor axis $g_{c,x}$ required to activate the wrinkling instability for various aspect ratios $a/b$ of the tube with fixed film thickness $w=0.02\,$mm and modulus ratio $\gamma^f/\gamma^s=10$. As seen in Fig.~\ref{fig:tu_psg}, there is a marginal increase in the critical growth $g_{c,x}$ with increasing aspect ratio. Comparing the critical buckling modes, it can be seen that increasing the aspect ratio leads to wrinkles not being formed in the regions where the substrate is thinner. This is attributed to the greater curvature of the film in regions of thinner substrate which makes it more resistant to wrinkling under circumferential compression when compared to a straighter film in regions where the substrate is thicker.    

\section{Summary}

We adopted a minimization-based variational framework for modeling the transient coupled problem of diffusion-driven finite elasticity, having the deformation map and fluid-volume flux as global primary fields. Upon spatial discretization of the incremental minimization problem, a variational-based structural stability analysis was carried out to judge the stability of a given chemo-mechanical state based on an eigenvalue analysis of the coupled global finite-element stiffness matrix. This procedure was implemented to investigate the onset of sinusoidal wrinkles on the surfaces of flat and tubular film-substrate hydrogel systems. The wrinkling instability mode was activated by subjecting the bilayers to swelling under suitable geometrical constraints. The variations in the critical buckling characteristics of the representative hydrogel systems were studied over a broad range of geometries and material parameters. For flat bilayers under lateral constraints, the dependence of the buckling characteristics on the system geometry resembled that of plates under uniaxial lateral compression. The trends showing the variations in critical load and wrinkle count with respect to the film-substrate modulus ratio were found to be in accordance with analytical and experimental studies available in literature on equivalent hyperelastic systems.

\textbf{Acknowledgement.}
The financial support of the German Research Foundation (DFG) within the Cluster of Excellence EXC 2075 (390740016) at the University of Stuttgart is gratefully acknowledged.

\appendix

\section{Conforming Raviart--Thomas-type finite-element design}
\label{app:A}
\begin{figure}[b]
	\centering
	\scriptsize
	\psfrag{aa}[c][c]{\footnotesize a)}
	\psfrag{bb}[c][c]{\footnotesize b)}
	\psfrag{phi}[c][c]{$\Bvarphi^h$}
	\psfrag{H}[l][l]{$\mathbb{H}^h$}
	\psfrag{xi}[c][c]{$\Bxi$}
	\psfrag{ntil}[c][c]{$\widetilde{\Bn}_0$}
	\psfrag{htil}[c][c]{$\widetilde{\mathbb{H}}^h$}
	\psfrag{xhat}[c][c]{$\widehat{\BX}(\Bxi)$}
	\psfrag{P}[c][c]{$\calP$}
	\psfrag{x}[l][l]{$\BX$}
	\psfrag{h}[l][l]{$\mathbb{H}^h$}
	\psfrag{n0}[l][l]{$\Bn_0$}
	\includegraphics[width=0.9\textwidth]{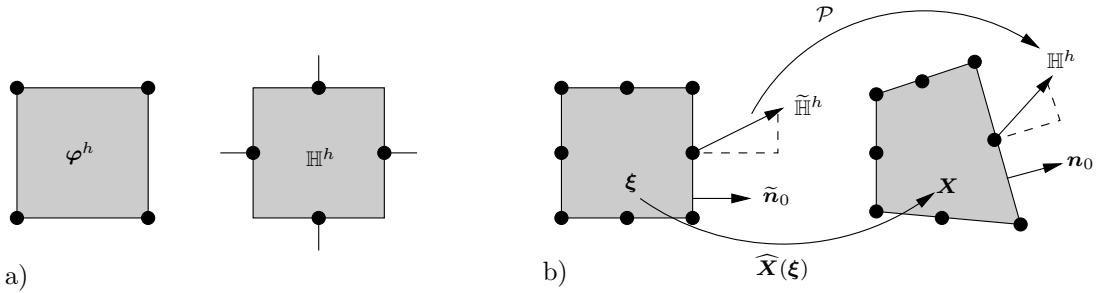}
	\captionsetup{width=0.95\textwidth}
	\caption{\footnotesize \textit{Structure of the Q\textsubscript{1}RT\textsubscript{0} element and the Piola transformation rule for $\mathbb{H}^h$.} a) The Q\textsubscript{1}RT\textsubscript{0} element consists of scalar-valued flux degrees of freedom, defined in (\ref{eq:2.28}), at the edge centers of a 4-noded Q\textsubscript{1} quadrilateral. These are interpolated by using piecewise linear vectorial RT\textsubscript{0} shape functions, defined in (\ref{eq:2.30}), to obtain the ansatz for the volume flux field. The corner nodes represent displacement degrees of freedom that are interpolated by the standard bilinear shape functions used for Q\textsubscript{1} elements. b) The Piola transformation rule $\calP$ is used to map the volume-flux field from the local parametric space to the physical space. It scales the normal trace of the volume flux as per (\ref{eq:2.31}), thereby ensuring that the values of scalar flux degrees of freedom remain unchanged upon transformation. The pictures are adopted from the work of \citet{teichtmeister2019}.}
	\label{fig:RT}
\end{figure}

The incremental minimization principle (\ref{eq:2.18}) demands the finite-element ansatz for the discretized volume-flux field $\mathbb{H}^h$ to lie in the space $\calH(\text{Div}, \calB_0)$. From the work of \citet{brezzifortin}, we infer that in order to ensure $\calH(\text{Div},\calB_0)$ conformity, the normal component $\mathbb{H}^h\cdot\Bn_0$ has to be continuous across each inter-element boundary in the global finite-element mesh. For two-dimensional problems, this continuity requirement is fulfilled by making use of the Q\textsubscript{1}RT\textsubscript{0} finite element for spatial discretization, which uses the lowest-order Raviart--Thomas-type vectorial interpolation functions for $\mathbb{H}^h$, first introduced in the work of \citet{raviart1977}. As seen in Fig.~\ref{fig:RT}a, the Q\textsubscript{1}RT\textsubscript{0} element constitutes the scalar quantity
\begin{equation}
	\label{eq:2.28}
	h^e_K=\int_{\calL^e_K}\mathbb{H}^h(\BX)\cdot\Bn_0 \,\mathrm{d}A\,,
\end{equation}     
representing the normal trace of the volume flux across edge $K$, as a degree of freedom at each edge center of a standard Q\textsubscript{1} quadrilateral. The subscript $0$ in RT\textsubscript{0} indicates that the quantity $\mathbb{H}^h\cdot\Bn_0$ is a polynomial of degree $0$, i.e.\ constant along each element edge. In the local parametric space of an element $e$, the ansatz for the volume flux $\widetilde{\mathbb{H}}^h(\Bxi)$ is then constructed by interpolating the nodal quantities $\widetilde{h}^e_K$ on the element edges by piecewise linear vectorial shape functions, i.e.
\begin{equation}
	\label{eq:2.29}
	\widetilde{\mathbb{H}}^h(\Bxi)=\sum_{K=1}^4 \mathrm{\textbf{N}}^K(\Bxi)\,\widetilde{h}^e_K\,,
\end{equation}  
where the vectorial shape functions $\mathrm{\textbf{N}}^K(\Bxi)$ take the form

\begin{equation}
	\label{eq:2.30}
	\mathrm{\textbf{N}}^1=
	\begin{bmatrix}
	0 \\[0.75ex] \frac{1}{4}(\xi_2+1)
	\end{bmatrix},\,\,
	\mathrm{\textbf{N}}^2=
	\begin{bmatrix}
	\frac{1}{4}(\xi_1+1) \\[0.75ex] 0
	\end{bmatrix},\,\,
	\mathrm{\textbf{N}}^3=
	\begin{bmatrix}
	0 \\[0.75ex] \frac{1}{4}(\xi_2-1)
	\end{bmatrix},\,\,
	\mathrm{\textbf{N}}^4=
	\begin{bmatrix}
	\frac{1}{4}(\xi_1-1) \\[0.75ex] 0
	\end{bmatrix}\,,
\end{equation}
refer \citet{teichtmeister2019}. Since the discrete scalar-valued fluxes represent degrees of freedom that are solved for in a finite-element context, their values should remain unchanged upon transformation from the parametric space to the physical space, i.e. $h^e_K \stackrel{!}{=} \widetilde{h}^e_K$. Since the normal trace of the flux field along an element edge is constant for RT\textsubscript{0}-type interpolation, this reduces to the condition
\begin{equation}
	\label{eq:2.31}
	\mathbb{H}^h(\BX)\cdot\Bn_0\,|\calL^e_K| \,\stackrel{!}{=}\,
	\widetilde{\mathbb{H}}^h(\Bxi)\cdot\widetilde{\Bn}_0\,|\widetilde{\calL}^e_K|\,,
\end{equation}
where $|\widetilde{\calL}^e_K|$ and $|\calL^e_K|$ represent the lengths of the edge $K$ in the parametric and physical spaces respectively. Therefore, the volume-flux field has to be transformed to the physical space in such a way that (\ref{eq:2.31}) is satisfied \textit{a priori}. This is achieved by adopting the Piola transformation rule $\calP$ (\cite{brezzifortin}, \cite{teichtmeister2019}) for the volume-flux field, which is expressed as
\begin{equation}
	\label{eq:2.32}
	\mathbb{H}^h(\BX)=\calP\left[ \widetilde{\mathbb{H}}^h(\Bxi)\right]  := \frac{1}{\widehat{J}(\Bxi)}\,\widehat{\BJ}(\Bxi)\left[ \widetilde{\mathbb{H}}^h(\Bxi)\right] \,.
\end{equation}   
Here, the matrix $\widehat{\BJ}(\Bxi):=\partial \widehat{\BX}/\partial \Bxi$ represents the Jacobian of the standard bilinear transformation $\widehat{\BX}(\Bxi)$ between parameteric and physical spaces for a Q\textsubscript{1} quadrilateral and $\widehat{J}(\Bxi)$ is its determinant. 

In order to ensure that volume outflux across a certain edge of an element in the global mesh equals volume influx in the neighbouring element sharing the same edge, a suitable sign convention needs to be adopted for the flux degrees of freedom on the edge centers. For this, we refer to the sign convention described and implemented in the works of \citet{anjam2015}, \citet{boeger2017} and \citet{teichtmeister2019}. Each Q\textsubscript{1}RT\textsubscript{0} element in the global mesh is traversed in the counter-clockwise direction. If upon doing so there is an increase in global node number of the corner nodes along a certain element edge, the flux degree of freedom on that edge for the element under consideration is assigned a positive sign. On the contrary, a decrease in the global corner-node number along an element edge results in the flux degree of freedom on that edge having a negative sign, which is taken into account by multiplying the vectorial RT\textsubscript{0} shape function corresponding to that edge by a factor of $-1$.

\end{document}